\input{amssymb.sty}

\documentclass[11pt]{amsart}
\theoremstyle{plain}
\newtheorem{thm}{Theorem}[section]
\newtheorem*{thm*}{Theorem}
\newtheorem{cor}[thm]{Corollary}
\newtheorem*{cor*}{Corollary}
\newtheorem*{conj*}{Conjecture}
\newtheorem*{lemma*}{Lemma}
\newtheorem{lemma}[thm]{Lemma}
\newtheorem*{prop*}{Proposition}
\newtheorem{prop}[thm]{Proposition}

\theoremstyle{definition}
\newtheorem{rems}[thm]{Remarks}
\newtheorem*{defn*}{Definition}
\newtheorem*{rems*}{Remarks}
\newtheorem*{proof*}{Proof}
\newtheorem{prel*}{Preliminaries}
\newtheorem{examples*}{Examples}

\newcommand{\C}{\mathbb{C}}
\newcommand{\E}{\widetilde{\mathcal E}}
\newcommand{\F}{\widetilde{\mathcal F}}
\newcommand{\M}{\widetilde{M}}
\newcommand{\W}{\widetilde{W}}
\newcommand{\J}{{\mathbb H}}

\newcommand{\npartial}{\not\partial}

\newcommand{\Ind}{\operatorname{Ind}}
\newcommand{\Index}{\operatorname{index}}

\newcommand{\rank}{\operatorname{rank}}
\newcommand{\tr}{\operatorname{tr}}
\newcommand{\Tr}{\operatorname{Tr}}

\newcommand{\End}{\operatorname{End}}

\newcommand{\R}{\operatorname{\mathbb R}}
\newcommand{\Z}{\operatorname{\mathbb Z}}

\newcommand{\nc}{\newcommand}
\nc{\nt}{\newtheorem}
\nc{\gf}[2]{\genfrac{}{}{0pt}{}{#1}{#2}}
\nc{\mb}[1]{{\mbox{$ #1 $}}}
\nc{\real}{{\mathbb R}}
\nc{\comp}{{\mathbb C}}
\nc{\ints}{{\mathbb Z}}
\nc{\Ltoo}{\mb{L^2({\mathbf H})}}
\nc{\rtoo}{\mb{{\mathbf R}^2}}
\nc{\slr}{{\mathbf {SL}}(2,\real)}
\nc{\slz}{{\mathbf {SL}}(2,\ints)}
\nc{\su}{{\mathbf {SU}}(1,1)}
\nc{\so}{{\mathbf {SO}}}
\nc{\hyp}{{\mathbb H}}
\nc{\disc}{{\mathbf D}}
\nc{\torus}{{\mathbb T}}

\nc{\ca}{{\mathcal A}}
\nc{\cag}{{{\mathcal A}^\Gamma}}
\nc{\cg}{{\mathcal G}}
\nc{\chh}{{\mathcal H}}
\nc{\ck}{{\mathcal B}}
\nc{\cl}{{\mathcal L}}
\nc{\cm}{{\mathcal M}}
\nc{\cs}{{\mathcal S}}
\nc{\cz}{{\mathcal Z}}
\nc{\sind}{\sigma{\rm -ind}}

\begin{document}

\title[TWISTED HIGHER INDEX THEORY ON GOOD ORBIFOLDS]
{TWISTED HIGHER INDEX THEORY ON GOOD ORBIFOLDS\\
AND FRACTIONAL QUANTUM NUMBERS}
\author{Matilde Marcolli}
\address{Department of Mathematics, Massachussetts Institute of Technology,
Cambridge, Mass., USA}
\email{matilde@math.mit.edu}
\author{Varghese Mathai}
\address{Department of Mathematics, University of Adelaide, Adelaide 5005,
Australia}
\email{vmathai@maths.adelaide.edu.au}


\begin{abstract}
In this paper, we study the twisted higher index theory
of elliptic operators
on orbifold covering spaces of compact good orbifolds,
which are invariant under a projective
action of the orbifold fundamental group, and we apply these results to obtain
qualitative results, related to generalizations of the
Bethe-Sommerfeld conjecture,
on the spectrum of  self adjoint elliptic operators
which are invariant under a projective action of the
orbifold fundamental group.
We also compute the
range of the higher traces on $K$-theory, which we then apply to
compute the
range of values of the Hall conductance  in the
quantum Hall effect on the hyperbolic plane. The new phenomenon
that we observe in this case is that the Hall conductance
again has plateaus at all energy levels belonging to
any gap in the spectrum of the Hamiltonian, where it
is now shown to be equal to an integral multiple of
a {\em fractional} valued
invariant. Moreover the set of possible denominators is finite
and has been explicitly determined. It is plausible that this
might shed light on the mathematical mechanism
responsible for fractional quantum numbers.
\end{abstract}

\maketitle

\section*{Introduction}
In this paper, we prove a twisted higher index theorem for
elliptic operators
on orbifold covering spaces of compact good orbifolds,
which are invariant under a projective
action of the orbifold fundamental group. These are basically
the evaluation of pairings of higher traces (which are
cyclic cocycles arising from the orbifold fundamental group
and the multiplier defining the projective action) with
the index of the elliptic operator, considered as an element
in the $K$-theory of some completion of the twisted
group algebra of the orbifold fundamental group.
The main purpose of generalizing the twisted higher index
theorem to orbifolds is to highlight the fact that
when the orbifold is not smooth, then the twisted higher
index can be a  {\em fraction}.
In the smooth case, the higher twisted index theorem was
used in \cite{CHMM} to study the quantum Hall effect
on hyperbolic space, and one of our key aims in this paper is to
generalize the results of \cite{CHMM} to general Fuchsian
groups and orbifolds, and can  be viewed equivalently
as the generalization
of results in \cite{CHMM} to the equivariant context.
As a result, we obtain a mathematical
mechanism that may explain the {\em fractional quantum
numbers}
that appear in the quantum Hall effect, cf. section 6.

Let $\Gamma\to\widetilde{M}\to M$ be a normal
covering space
of a compact smooth manifold $M$. Then in [A] (and clarified by [CM])
Atiyah showed
that any $\Gamma$-invariant elliptic differential operator $\widetilde{P}$
acting on
$L^2$ sections on $\widetilde{M}$, and which is the lift of an elliptic
differential
operator $P$ on $M$, yields via the parametrix construction, an element in
$K$-theory
\[
\Ind_\Gamma (\widetilde{P})\in K_0(\mathcal{R}(\Gamma)) \] where
$\mathcal{R}(\Gamma)$ is the group ring of $\Gamma$ with coefficients in
the algebra
of rapidly decreasing matrices, i.e. \[
\mathcal{R}=\left\{\left(a_{ij}\right)_{i,j\in\mathbb{N}}\,:\,
\sup_{i,j\in\mathbb{N}} i^k j^\ell \left|a_{ij}\right|<\infty
\ \forall\ k,\ell\in\mathbb{N}\right\}
\] Let $\tr$ denote the trace on $\mathbb{C}$ the group algebra
$\mathbb{C}(\Gamma)$. . More precisely,
if $\delta_g, \, \, g\in \Gamma$ denotes the canonical basis of
$\mathbb{C}$, i.e. $\delta_g(g') = 1$ if $g = g'$ and zero otherwise, then
\[
\tr\left(\delta_g\right) = \begin{cases} 1 & \text{ if } g=e \\  0 & \text{
otherwise } \end{cases}
\] and $\Tr\,:\,\mathcal{R}\to\mathbb{C}$ denote the trace on $\mathcal{R}$,
$\Tr\left(\left(a_{ij}\right)_{i,j\in\mathbb{N}}\right)
=\sum_{i\in\mathbb{N}} a_{ii}$.

Observe that the natural inclusion map $j: \mathcal{R}(\Gamma) \to
C^*(\Gamma)$ induces a morphism in $K$-theory
$$
j_*: K_\bullet (\mathcal{R}(\Gamma) ) \to K_\bullet (C^*(\Gamma))
$$

Then the cup product \ $\tr\sharp\Tr$\  is a trace on $\mathcal{R}(\Gamma)$
which
extends to $K$-theory. We also denote by  \ $\tr\sharp\Tr$\  the canonical
trace on
$C^*(\Gamma)$. Atiyah \cite{At} then proved that \[
[\tr\sharp\Tr]\left(\Ind_\Gamma(\widetilde{P})\right) =
[\tr\sharp\Tr]\left(j_*(\Ind_\Gamma(\widetilde{P}))\right) = \Index(P) \] where
$\Index(P)$ denotes the Fredholm index of the elliptic operator $P$. Using the
Atiyah-Singer index theorem, he was then able to establish a cohomological
formula for
$[\tr\sharp\Tr]\left(\Ind_\Gamma(\widetilde{P})\right)$.

Since then, there have been two significant generalizations of this
theorem. The
first is due to Connes and Moscovici \cite{CM}, where they compute the
pairing of
$\Ind_\Gamma(\widetilde{P})$ with other higher traces. More precisely,
given a
normalized group cocycle $c\in Z^k(\Gamma,\mathbb{C})$, they define a
cyclic cocycle
$\tr_c$ on the group ring $\mathbb{C}(\Gamma)$ as follows:
\[
\tr_c\left(\delta_{g_0},\dots,\delta_{g_k}\right) = \begin{cases}
c\left(g_1,\dots,g_k\right) & \text{ if } g_o  g_1,\dots,g_k = 1 \\ 0 & \text{
otherwise }\end{cases}
\] Then the cup product \ $\tr_c\sharp\Tr$ \ extends continuously to a
k-dimensional cyclic
cocycle on
$\mathcal{R}(\Gamma)$ which extends to $K$-theory, and \cite{CM} establish
a cohomological formula for \
$[\tr_c\sharp\Tr]\left(\Ind_\Gamma(\widetilde{P})\right)$.

Now let $\sigma$ be a multiplier on $\Gamma$ and suppose that there is a
projective
$(\Gamma,\bar\sigma)$ action on $L^2$ sections on $\widetilde{M}$, cf
section 1.3. Then in
\cite{Gr}, Gromov
extends Atiyah's index theorem in another direction, to an index theorem
for elliptic
operators ${D}$ on $\widetilde{M}$ which are invariant under the projective
$(\Gamma,\bar\sigma)$ action. More precisely, Gromov essentially remarked
(and clarified in this paper) that one could
modify Atiyah's parametrix construction to obtain a
$(\Gamma, \sigma)$-index which is an element in $K$-theory \[
\Ind_{(\Gamma,\sigma)} ({D})\in K_0(\mathcal{R}(\Gamma,\sigma)) \] where
$\mathcal{R}(\Gamma,\sigma)$ is the group ring of $\Gamma$ with coefficients in
$\mathcal{R}$, and which is twisted by the multiplier $\sigma$. Let $\tr$
denote the trace on $\mathbb{C}(\Gamma,\sigma)$. More precisely,
if $\delta_g, \, \, g\in \Gamma$ denotes the canonical basis of
$\mathbb{C}$, i.e. $\delta_g(g') = 1$ if $g = g'$ and zero otherwise, then
\[
\tr(\delta_g)=\begin{cases} 1 & \text{ if } g=e \\  0 & \text{ otherwise
}.\end{cases}
\] Then the cup product \ $\tr\sharp\Tr$\  is a trace on
$\mathcal{R}(\Gamma,\sigma)$
which extends to $K$-theory. Observe that the natural inclusion map $j:
\mathcal{R}(\Gamma, \sigma) \to
C^*(\Gamma, \sigma)$ induces a morphism in $K$-theory
$$
j_*: K_\bullet (\mathcal{R}(\Gamma,\sigma) ) \to K_\bullet (C^*(\Gamma,\sigma))
$$
Gromov also computes a cohomological formula for the $(\Gamma, \sigma)$-index
\[ [\tr\sharp\Tr]\left(\Ind_{(\Gamma,\sigma)}({D})\right)
= [\tr\sharp\Tr]\left(j_*(\Ind_{(\Gamma,\sigma)}({D}))\right). \]

In this paper,
we will prove an index theorem which will generalize and unify the index
theorems of
Atiyah, Connes and Moscovici, and Gromov, in the case of
good orbifolds, that is orbifolds such that their
orbifold universal cover is a
manifold.
Now let $\Gamma\to\widetilde{M}\to M$ denote the universal orbifold
cover of a compact good orbifold $M$, so that $\widetilde M$ is a
smooth manifold.
Let $\sigma$ be a multiplier on
$\Gamma$ and assume that there is a projective $(\Gamma,\bar\sigma)$-action
on $L^2$ sections of $\Gamma$-invariant vector bundles over $\widetilde M$.
By considering $(\Gamma,\bar\sigma)$-invariant elliptic operators acting on
$L^2$ sections of these bundles, we will show that again, this defines
$(\Gamma, \sigma)$-index element in $K$-theory
\[
\Ind_{(\Gamma,\sigma)} ({D})\in K_0(\mathcal{R}(\Gamma,\sigma)). \]
We will compute the pairing of
$\Ind_{(\Gamma,\sigma)}({D})$ with higher traces. More precisely,
given a normalized group cocycle $c\in Z^k(\Gamma, \mathbb C)$, $
k=0,\ldots, \dim M$,
we define a cyclic cocycle
$\tr_c$ of dimension $k$ on the twisted group ring
$\mathbb C(\Gamma, \sigma)$, which is
given by
$$
\tr_c(a_0\delta_{g_0}, \ldots,a_k\delta_{g_k}) = \left\{\begin{array}{l}
a_0\ldots a_k c(g_1,\ldots,g_k)
\tr( \delta_{g_0}\delta_{g_1} \ldots\delta_{g_k}) \quad\text{if} \,\, g_0
\ldots g_k =1\\
\\
0 \qquad \text{otherwise.}
\end{array}\right.
$$
where $a_j \in \C$ for $j=0,1,\ldots k$.
Of particular interest is the case when $k=2$, when the formula above
reduces to
$$
\tr_c(a_0\delta_{g_0},a_1\delta_{g_1} ,a_2\delta_{g_2}) =
\left\{\begin{array}{l}
a_0 a_1 a_2 c(g_1, g_2)
\sigma(g_1, g_2) \quad \text{if} \,\, g_0 g_1 g_2 =1 ;\\
\\
0 \qquad \text{otherwise.}
\end{array}\right.
$$
The cup product $\tr_c \#\Tr$ extends
continuously to a $k$-dimensional cyclic cocycle on $\mathcal{R}(\Gamma,
\sigma)$, which
then extends
to $K$-theory. We will also compute a cohomological formula for
$$
[\tr_c\sharp\Tr]\left(\Ind_{(\Gamma,\sigma)}({D})\right).
$$
Our method consists of applying the Connes-Moscovici
{\em local} higher index theorem to a
family of idempotents constructed from the heat operator on $\widetilde{M}$,
all of which represent the $(\Gamma, \sigma)$-index.
It is interesting to mention that the orbifold case differs from the
smooth case: the index and the $L^2$-index are different due to the
presence of the singular stratum. The contribution of the singular
stratum is present in the index formula, but is not detected by the
$L^2$-index.

Let $\Gamma$ be a Fuchsian group of signature $(g, \nu_1,\ldots,
\nu_n)$ (cf. section 1 for more details),
that is, $\Gamma$ is the orbifold fundamental group of
the 2 dimensional hyperbolic orbifold $\Sigma(g, \nu_1,\ldots,\nu_n)$
of signature $(g, \nu_1,\ldots,\nu_n)$.
Using a result of Kasparov \cite{Kas1} on $K$-amenable groups
as well as a calculation by Farsi \cite{Far} of the orbifold $K$-theory
$K^\bullet_{orb} (\Sigma(g, \nu_1,\ldots,\nu_n))$ of compact 2-dimensional
hyperbolic orbifolds, we are able to compute
the $K$-theory of twisted group $C^*$ algebras, under the
assumption that the Dixmier-Douady invariant of the multiplier
$\sigma$ is trivial
$$
K_j(C^*(\Gamma, \sigma)) \cong \left\{\begin{array}{l}
\Z^{2-n +\sum_{j=1}^n \nu_j}\qquad \text{if}\,\, j=0;\\ \\
\mathbb{Z}^{2g} \qquad \qquad\qquad\text{if}\,\, j=1.
\end{array}\right.
$$
Notice that  $K_0$ is much larger in the general Fuchsian group
case than in the
torsionfree case, where $K_0$ was determined to be always
$\Z^2$,  \cite{CHMM}.
We also show that the orbifold $K$-theory of
any 2-dimensional orbifold is
generated by orbifold line bundles.
The result is derived by means of equivariant $K$-theory and the Baum-Connes
\cite{BC} equivariant Chern character with values in the
delocalized equivariant cohomology of the smooth surface $\Sigma_{g'}$
that covers the good orbifold $\Sigma(g, \nu_1,\ldots,\nu_n)$. We show
that the Seifert invariants, cf. \cite{Sc}, correspond to
the pairing of the equivariant Chern character of \cite{BC}
with a fundamental class in the delocalized equivariant homology of
$\Sigma_{g'}$.

Using these results and our twisted higher index theorem for orbifolds, we
compute in section 3 under the same assumptions as before,
the range of the trace on $K$-theory to be,
$$
[\tr\sharp\Tr]( K_0 (C^*(\Gamma,\sigma)))
= \Z\theta + \Z + \sum_{i=1}^n \Z (1/\nu_i)
$$
where $\theta$ denotes the evaluation of the multiplier $\sigma$ on the
fundamental class of $\Gamma$.
We then apply our calculation of the range of the trace on $K$-theory
to study some quantitative aspects of the spectrum of projectively
periodic elliptic operators on the hyperbolic plane.
Some of the most outstanding open problems about magnetic Schr\"odinger
operators or Hamiltonians
on Euclidean space is concerned with the nature of their spectrum,
and and are the {\em Bethe-Sommerfeld conjecture} (BSC) and the
{\em Ten Martini Problem} (TMP) (cf. \cite{Sh}). More
precisely,
TMP asks whether
given a multiplier $\sigma$ on $\mathbb Z^2$, is there an associated
Hamiltonian
({\it i.e.} a Hamiltonian which commutes with the $({\mathbb Z}^2,
\bar\sigma)$
projective
action of $\mathbb Z^2$ on $L^2(\mathbb R^2)$) possessing a Cantor set
type spectrum, in the sense that the intersection of
the spectrum of the Hamiltonian with some compact interval in $\mathbb R$
is a Cantor set?
One can deduce from the range of the
trace on $K_0$ of the twisted group $C^*$-algebras that when the
multiplier takes its values in the  roots of unity in $U(1)$ (we say then that
it is rational) that such a Hamiltonian cannot exist.
However, in the Euclidean
case and for almost all irrational numbers, the discrete form of TMP has
been settled
in the affirmative, cf.
\cite{Last}.
BSC asserts that if the multiplier is trivial, then the spectrum
of any associated Hamiltonian has only a {\em finite} number of
gaps. This was first established in the Euclidean case by Skrigonov
\cite{Skri}.
In Sections 3 and 4, we are concerned also with generalizations
of the TMP and the BSC, which we call the
{\em Generalized Ten Dry Martini Problem}  and the
{\em Generalized Bethe-Sommerfeld
conjecture}. We
prove that the Kadison constant
 of the twisted group $C^*$-algebra $C^*_r(\Gamma,\sigma)$
is positive whenever the multiplier is rational, where $\Gamma$ is now
the orbifold fundamental
group of a signature $(g, \nu_1, \ldots ,\nu_n)$ hyperbolic orbifold.
We then use the results of Br\"uning and Sunada \cite{BrSu}
to deduce that when the multiplier is rational, the
generalized Ten Dry Martini Problem is answered in the
negative, and we leave open the more difficult irrational case.
More precisely, we show that the spectrum of such
a $(\Gamma, \bar\sigma)$ projectively periodic elliptic
operator is the union of countably many (possibly degenerate)
closed intervals which can only accumulate at infinity.
This also gives evidence that the generalized
Bethe-Sommerfeld conjecture is true,
and generalizes earlier results of \cite{CHMM} in the torsion-free case.
In section 4, we again use the range of the trace theorem above, together
with other geometric arguments to give a complete classification
upto isomorphism of the twisted group $C^*$ algebras
$C^*(\Gamma,\sigma)$, where $\sigma$ is assumed to have
trivial Dixmier-Douady invariant as before.

In section 5, we use a result of \cite{Ji},
which is a twisted analogue of a result of Jollissant
and which says in particular that when $\Gamma$ is a cocompact Fuchsian
group, then
the natural inclusion map $j: \mathcal{R}(\Gamma,\sigma) \to
C^*(\Gamma,\sigma)$ induces an isomorphism in $K$-theory
$$
K_\bullet (\mathcal{R}(\Gamma,\sigma) ) \cong K_\bullet (C^*(\Gamma,\sigma))
$$
Using this, together with our twisted higher index theorem for good
orbifolds, and under the same assumptions as before,
we are able to compute the range of the higher trace on $K$-theory
$$
[\tr_c\sharp\Tr]( K_0 (C^*(\Gamma,\sigma))) = \phi\Z
$$
where $\phi  =
2(g-1)+ (n-\nu) \in {\mathbb Q}, \quad \nu = \sum_{j=1}^n 1/\nu_j$ and
$c$ is the area 2-cocycle on $\Gamma$, i.e. $c$ is the restriction to $\Gamma$
of the area 2-cocycle on $PSL(2, \R)$, cf. section 5.
We will give examples of good
2-dimensional orbifolds in section 5 for which
$\phi$ is a {\em fraction}; however it is an
{\em integer} whenever the orbifold is smooth, i.e. whenever $1=\nu_1=\ldots
=\nu_n$, which was the case that was considered in \cite{CHMM}.

In section 6, we study
the hyperbolic Connes-Kubo formula for the
Hall conductance in
the discrete model of the Quantum Hall Effect on the hyperbolic plane,
where we consider Cayley graphs of Fuchsian groups which may have
torsion subgroups, generalizing results in \cite{CHMM} where
only torsion-free Fuchsian groups were considered. We recall that
the results in \cite{CHMM} generalised to
hyperbolic space the noncommutative geometry
approach to the Euclidean
quantum Hall effect that was pioneered by
Bellissard and collaborators \cite{Bel+E+S}, Connes \cite{Co}
and Xia \cite{Xia}.
The Cayley graphs of these Fuchsian groups
are not in  general trees, as they may now have loops.
We first relate the hyperbolic Connes-Kubo cyclic 2-cocycle
and the area cyclic 2-cocycle on the algebra ${\mathcal R} (\Gamma, \sigma)$,
and show that they define the same class in cyclic cohomology.
Then we use use the range of the higher trace on $K$-theory to determine
the range of values of the Hall conductance in the Quantum Hall Effect.
The new phenomenon
that we observe in this case is that the Hall conductance
again has plateaus at all energy levels belonging to
any gap in the spectrum of the Hamiltonian (known as the generalized
Harper operator), where it
is now shown to be equal to an integral multiple of
a {\em fractional} valued
invariant. Moreover the set of possible denominators is finite
and has been explicitly determined, \cite{Bro}. It is plausible that this
might shed light on the mathematical mechanism
responsible for fractional quantum numbers.

\noindent{\bf Acknowledgments:} The second author thanks A. Carey and K.
Hannabus
for some clarifying comments concerning section 6.

\section{Preliminaries}

\subsection{Good orbifolds}
For further details on the fundamental material
on orbifolds, see \cite{Sc}, \cite{FuSt} and \cite{Bro}.

The definition of an orbifold generalises that of a manifold.
More precisely,
an {\em orbifold} $M$ of dimension $m$ is a Hausdorff, second countable
topological space with a Satake atlas $\mathcal V = \{U_i, \phi_i\}$
which covers $M$, consisting of open sets $U_i$ and homeomorphisms $\phi_i
: U_i
\to D^m/G_i$, where $D^m$ denotes the unit ball in ${\R}^m$
and $G_i$ is a finite subgroup of the orthogonal group $O(m)$,
satisfying the following compatibility relations; the compositions
$$
\phi_j\circ\phi_i^{-1} : \phi_i(U_i\cap U_j) \to \phi_j(U_i\cap U_j)
$$
locally lifts to be a smooth map $\R^m \to \R^m$, whenever
the intersection $U_i\cap U_j \ne \emptyset$.
The open
sets $U_i$ are called local orbifold charts.
In general, an orbifold $M$ can be obtained as a quotient
$M=X/G$ of an infinitesimally free compact Lie group action on a smooth
manifold $X$. In fact, by Satake \cite{Sat} and Kawasaki \cite{Kaw},
$X$ can chosen to be the smooth manifold of orthonormal frames
of the orbifold tangent bundle of $M$ (cf. section 1.4) and
$G$ can be chosen to be the orthogonal group $O(m)$.

An orbifold covering of $M$ is an orbifold map $f: Y\to
M$, where $Y$ is also an orbifold, such that any point on $M$ has a
neighbourhood $U$ such
that $f^{-1}(U)$ is the disjoint union of open sets $U_\alpha$, with
$f\mid_{U_\alpha}: U_\alpha \to U$  a quotient map between two
quotients of $\R^k$ by finite groups $H_1 < H_2$. The generic
fibers of the covering map $f$ are isomorphic to a discrete
group which acts as deck transformations.

An orbifold  $M$ is {\em good} if it is
orbifold-covered by a smooth manifold; it is {\em bad} otherwise.
A good orbifold is said to be {\em orientable} if it is
orbifold covered
by an oriented manifold and the deck transformations act via
orientation preserving diffeomorphisms on the orbifold cover.
Equivalently, as shown in \cite{Sat} and \cite{Kaw},
an orbifold is orientable if it has an oriented frame bundle
$X$ such that $M = X/ SO(m)$.

We next recall briefly some basic
notions on Euclidean and hyperbolic orbifolds, which are
by fiat orbifolds whose universal orbifold covering space is
Euclidean space and hyperbolic space respectively. We are mainly
interested in the case of 2 dimensions, and
we will assume that the orbifolds in this paper are orientable.

A 2-dimensional compact orbifold has singularities that are cone
points or reflector lines. Up to passing to $\Z_2$-orbifold covers, it
is always possible to reduce to the case with only isolated cone
points.

Let $\mathbb H$ denote the hyperbolic
plane and $\Gamma$ a Fuchsian group of signature $(g, \nu_1,\ldots,
\nu_n)$, that is, $\Gamma$ is a discrete cocompact subgroup of
$PSL(2, \R)$ of genus $g$ and with $n$ elliptic elements of order
$ \nu_1,\ldots, \nu_n$ respectively. Explicitly,
$$
\Gamma = \Big\{A_i, B_i, C_j \in PSL(2, \R)| i=1,\ldots g, \quad
j = 1, \ldots n, $$
$$
\prod_{i=1}^g[A_i, B_i] C_1\ldots C_n = 1, \quad C_j^{\nu_j} =1, \quad
j = 1, \ldots n\big\}
$$
Then the corresponding compact oriented hyperbolic 2-orbifold of signature
$(g, \nu_1,\ldots, \nu_n)$ is defined as the quotient space
$$
\Sigma(g, \nu_1,\ldots,\nu_n)= \Gamma\backslash {\mathbb H}.
$$
A compact oriented 2-dimensional Euclidean orbifold is obtained
in a similar manner, but with $\mathbb H$ replaced by $\R^2$, and
a complete list of these can be found in \cite{Sc}.
Then $\Sigma(g, \nu_1,\ldots,\nu_n)$ is a compact surface of
genus $g$ with $n$
elliptic points $\{p_j\}_{j=1}^n$ such that each $p_j$ has a small
coordinate neighborhood $U_{p_j} \cong D^2/{\Z_{\nu_j}}$, where $D^2$
denotes the unit disk in $\R^2$ and $\Z_{\nu_j}$ is the cyclic
group of order ${\nu_j}, \quad j = 1, \ldots n$. Observe that
the complement $\Sigma(g, \nu_1,\ldots,\nu_n) \setminus
\cup_{j=1}^n U_{p_j}$ is a compact Riemann surface of genus
$g$ and with $n$ boundary components. The group $\Gamma$ is the
orbifold fundamental group of $\Sigma(g, \nu_1,\ldots,\nu_n)$, where
the generators $C_j$ can be represented by the $n$ boundary components
of the surface $\Sigma(g, \nu_1,\ldots,\nu_n) \setminus
\cup_{j=1}^n U_{p_j}$.

All Euclidean and hyperbolic 2-dimensional
orbifolds $\Sigma(g,\nu_1,\ldots,\nu_n)$ are
good, being in fact orbifold covered by a smooth surface
$\Sigma_{g'}$ cf. \cite{Sc}, i.e. there is a finite group $G$ acting on
$\Sigma_{g'}$ with quotient $\Sigma(g,\nu_1,\ldots,\nu_n)$,
where $\#(G) = 1+ \sum_{j=1}^n (\nu_j-1)$ and $g' =
1+ \frac{\#(G)}{2}(2(g -1) + (n - \nu))$ and where
$\nu = \sum_{j=1}^n 1/\nu_j$.
According to the classification of 2-dimensional orbifolds
given in \cite{Sc}, the only bad 2-orbifolds are the ``teardrop'',
with underlying surface $S^2$ and one cone point of cone angle
$2\pi/p$, and the ``double teardrop'', with underlying surface $S^2$
and two cone points with angles $2\pi/p$ and $2\pi/q$, $p\neq q$.

In this paper we restrict our attention to good orbifolds. It should
be pointed out that the techniques used in this paper cannot be
extended directly to the case
of bad orbifolds. It is reasonable to expect that the index theory on
bad orbifolds will involve analytical techniques for more general
conic type singularities.

\subsection{On Groupoids}
Recall that a groupoid consists of a set $\mathcal{G}$ together with a
distinguished subset $\mathcal{G}^{(0)}\subset\mathcal{G}$ and two maps
\[
   r,s : \mathcal{G}\to\mathcal{G}^{(0)}
\]
and a composition law
\[
   \circ : \mathcal{G}^{(2)}=\left\{(\gamma_1,\gamma_2)\in\mathcal{G}
    \times\mathcal{G} : s(\gamma_1)=r(\gamma_2)\right\}\to\mathcal{G}
\]
such that
\begin{enumerate}
\item $s(\gamma_1\circ\gamma_2)=s(\gamma_2)$ and $r(\gamma_1\circ\gamma_2)=
\gamma_1)\forall\gamma_1,\gamma_2\in\mathcal{G}^{(2)}$
\item $s(x)=r(x)=x\quad\forall x\in\mathcal{G}^{(0)}$
\item $\gamma\circ s(\gamma)=\gamma$ and $r(\gamma)\circ\gamma=\gamma\quad
\forall\quad\gamma\in\mathcal{G}$
\item $(\gamma_1\circ\gamma_2)\circ\gamma_3=\gamma_1\circ(\gamma_2\circ
\gamma_3)$
\item Each $\gamma$ has a 2-sided inverse $\gamma^{-1}$ such that
$\gamma\circ\gamma^{-1}=r(\gamma)$ and $r^{-1}\circ \gamma=
\nolinebreak[4]s(\gamma)$
\end{enumerate}
$r$ is usually called the {\em range map} and $s$ the {\em source map}
of the groupoid $\mathcal{G}$.

\begin{enumerate}

\item $\mathcal{G}=\text{group},\ \mathcal{G}^{(0)}=\{e\}$ and the range
and source maps are degenerate;
\item The orbifold fundamental groupoid of $M$. Let $M$ be a good
orbifold, $\tilde M$ be the orbifold universal cover, and $\Gamma$ be
the orbifold fundamental group. We have
\begin{gather*}
   \mathcal{G} = \widetilde{M} \times_{\Gamma} \widetilde{M},\quad
     \mathcal{G}^{(0)} = \widetilde{M}/\Gamma = M \\
   r(\tilde{x},\tilde{y})=x\in M,\quad s(\tilde{x},\tilde{y})=y\in M \\
   (\tilde{x},\tilde{y})\circ(\tilde{y},\tilde{z})=(\tilde{x},\tilde{z}).
\end{gather*}

\end{enumerate}

\noindent Let $u\in\mathcal{G}^{(0)}$ and define $\mathcal{G}^u=r^{-1}(u),
\mathcal{G}_u=s^{-1}(u)$
\[
   \mathcal{G}^u_u=\{\gamma\in\mathcal{G}:u=r(\gamma)=s(\gamma)\}
     =\mathcal{G}^u\cap\mathcal{G}_u
\]
in the example 2. above, $\mathcal{G}^u=\widetilde{M}=\mathcal{G}_u$ and
$\mathcal{G}^u_u=\Gamma$.

The orbifold fundamental groupoid of a good orbifold is
the main example that we will be concerned with in this paper.

\subsection{Twisted Groupoid $C^*$ algebras}
Let $M$ be a good, compact orbifold, and $\mathcal E \to M$ be an orbifold
vector bundle over $M$, and $\E \to \widetilde{M}$ be its lift to
the universal orbifold covering space $\Gamma\to \widetilde{M}\to M$,
which is by assumption
a simply-connected smooth manifold.
We will now briefly review how to construct a
$(\Gamma, \bar\sigma)$-action (where $\sigma$ is a multiplier on
$\Gamma$ and $\bar\sigma$ denotes its complex conjugate)
on $L^2(\M)$. Let $\omega = d\eta$ be an exact $2$-form on $\widetilde{M}$
such that $\omega$ is also $\Gamma$-invariant, although $\eta$ is {\em not}
assumed to be $\Gamma$-invariant.
Define a Hermitian connection on the trivial line bundle over $\widetilde{M}$
as
\[
   \nabla = d + i\eta
\]
Its curvature is $\nabla^2=i\omega$.  Then $\nabla$ defines a
$(\Gamma,\bar{\sigma})$ action on $E$
$L^2(\widetilde{M},\widetilde{\mathcal{S}^\pm\otimes E})$ as follows:

Since $\omega$ is $\Gamma$ invariant, one has $\forall\gamma\in\Gamma$
\[
   0 = \gamma^* \omega-\omega = d(\gamma^* \eta-\eta),
\]
so that $\gamma^*\eta-\eta$ is a closed $1$-form on a simply-connected manifold
$\widetilde{M}\Rightarrow \gamma^*\eta-\eta=d\psi_\gamma$ for some smooth
function $\psi_\gamma$ on $\widetilde{M}$ satisfying
\begin{align*}
   \bullet\quad & \psi_\gamma(x)+\psi_{\gamma'}(\gamma
     x)-\psi_{\gamma'\gamma}(x) \qquad\text{is independent of $x\ \forall
       x\in\widetilde{M}$}\\
   \bullet\quad & \psi_\gamma(x_0)=0\qquad\text{for some $x_0\in\widetilde{M}
     \ \ \forall \gamma\in\Gamma$}
\end{align*}
Then
$\bar{\sigma}(\gamma,\gamma')=\exp(i\psi_\gamma(\gamma' x_0))$
defines a multiplier on $\Gamma$.
 Now define the
$(\Gamma,\bar{\sigma})$ action as follows:
For $u\in L^2(\widetilde{M}, \widetilde{\mathcal{S}^\pm\otimes E}),
      \quad U_\gamma\ u=\gamma^* u, \quad
S_\gamma u = \exp(i\,\psi_\gamma)\ u$, define $T_\gamma=U_\gamma\circ
S_\gamma$.
Then it satisfies $T_{\gamma_1} T_{\gamma_2}
=\bar{\sigma}(\gamma_1,\gamma_2)\,T_{\gamma_1 \gamma_2}$, and so it defines
a $(\Gamma,\bar{\sigma})$-action. It can be shown that
only multipliers $\bar \sigma$ such that the Dixmier-Douady invariant
$ \delta(\bar\sigma) = 0$ can give
rise to $(\Gamma, \bar\sigma)$-actions in this way cf. section 3.2
for a further discussion.

Let $D: L^2(\widetilde{M}, \E) \to L^2(\widetilde{M}, \E)$ be a self
adjoint elliptic differential
operator that commutes with a $(\Gamma, \bar\sigma)$-action
$T_\gamma \quad \forall \gamma\in \Gamma$ on $L^2(\widetilde{M}, \E)$.
Then by the functional calculus, all the spectral projections of $D$,
$E_\lambda = \chi_{[0,\lambda]}(D)$ are bounded self adjoint operators
on $L^2({\widetilde M}, \E)$ that commute with $T_\gamma \quad \forall
\gamma\in \Gamma$.
Now the commutant of the $(\Gamma, \bar\sigma)$-action on $L^2({\widetilde
M}, \E)$
is a von Neumann algebra
$$
{\mathcal U}(\mathcal{G}, \sigma) = \left\{Q\in B(L^2({\widetilde M}, \E))|
\quad T_\gamma Q =
Q T_\gamma \quad \forall \gamma\in \Gamma \right\}
$$
where $\mathcal{G}$ denotes the orbifold homotopy groupoid of
the orbifold $M$.
Since $T_\gamma Q = Q T_\gamma $, one sees that
$$
e^{i\phi_\gamma(x)}k_Q(\gamma x,\gamma y) e^{-i\phi_\gamma(y)} = k_Q(x,y)
$$
$\forall x, y \in {\widetilde M}\quad \forall \gamma\in \Gamma$,
where $k_Q$ denotes the Schwartz kernel of $Q$. In particular, observe
that $\tr(k_Q(x,x))$ is a $\Gamma$-invariant function on ${\widetilde M}$.
Using this, one sees that there is a semifinite trace on this
von Neumann algebra
$$
\tr: {\mathcal U}(\mathcal{G}, \sigma) \to \C
$$
defined as in the untwisted case due  to Atiyah \cite{At},
$$
Q \to \int_{M} \tr(k_Q(x,x))dx
$$
where $k_Q$ denotes the Schwartz kernel of $Q$. Note that this trace is finite
whenever $k_Q$ is continuous in a neighborhood of the diagonal in
${\widetilde M}\times {\widetilde M}$.
We remark that ${\mathcal U}(\mathcal{G}, \sigma)$ can also be defined as
the weak $*$-completion of a twisted convolution algebra of functions
on the groupoid $\mathcal{G}$, but we will not have use for this
alternate description here.

By elliptic regularity, the spectral projection
$E_\lambda$ has a smooth Schwartz kernel, so that
in particular, the spectral density function,
$N_Q(\lambda) = \tr(E_\lambda) <\infty \quad \forall \lambda$, is well defined.

If $\mathcal F$ is a fundamental domain for the action of $\Gamma$ on
${\widetilde M}$,
one sees that
$$
L^2({\widetilde M}, \E) \cong L^2(\Gamma)\otimes L^2(\mathcal F,
\E|_{\mathcal F})
$$
which can be proved by choosing a bounded measurable almost everywhere smooth
section of the orbifold covering ${\widetilde M} \to M$. Therefore it
follows easily that
$$
{\mathcal U} (\mathcal{G}, \sigma) \cong
{\mathcal U}(\Gamma, \sigma)\otimes  B(L^2(\mathcal F, \E|_{\mathcal F}))
$$
where the twisted group von Neumann algebra ${\mathcal U}(\Gamma, \sigma)$
is the weak closure of the twisted group algebra
$\C(\Gamma, \sigma)$ and $B(L^2(\mathcal F, \E|_{\mathcal F}))$
denotes the algebra of bounded operators on the Hilbert space
$L^2(\mathcal F, \E|_{\mathcal F})$.
There is a natural subalgebra ${C^*} (\mathcal{G}, \sigma)$
of ${\mathcal U} (\mathcal{G}, \sigma)$ which is defined
as follows. Let
\begin{align*}
{C_c}^\infty (\mathcal{G}, \sigma)=
& \{Q\in {\mathcal U} (\mathcal{G}, \sigma)
| k_Q \ \text{is smooth and supported in a compact }\\
&\text{neighborhood
of the diagonal}\}
\end{align*}
Then ${C^*}(\mathcal{G}, \sigma)$ is defined to be the {\em norm closure}
of ${C_c}^\infty (\mathcal{G}, \sigma)$.
It can also be shown to be the norm closure of
$$
\left\{Q\in {\mathcal U}_{{\widetilde M}} (\Gamma, \sigma)
| k_Q \ {\text{is smooth and}}\  k_Q(x,y) \ {\text {is $L^1$ in both the
$x$ and $y$
variables seperately}}\right\}
$$
The elements of ${C^*} (\mathcal{G}, \sigma)$
have the additional property of
some off-diagonal decay. In the earlier notation,
it can be shown that (cf. \cite{MRW})
$$
{C^*} (\mathcal{G}, \sigma) \cong {C^*}(\Gamma, \sigma)
\otimes {\mathcal K}(L^2(\mathcal F, \E|_{\mathcal F}))
$$
where the twisted group $C^*$ algebra ${C^*}(\Gamma, \sigma)$
is the norm closure of the twisted group algebra
$\C(\Gamma, \sigma)$ and ${\mathcal K}(L^2(\mathcal F, \E|_{\mathcal F}))$
denotes the algebra of compact operators on the Hilbert space
$L^2(\mathcal F, \E|_{\mathcal F})$.
We remark that ${C^*}(\mathcal{G}, \sigma)$ can also be defined as
the  norm completion of a twisted convolution algebra of functions
on the groupoid $\mathcal{G}$, but we will not have use for this
alternate description here.

\subsection{The $C^*$ algebra of an orbifold}
Let $M$ be an oriented orbifold of dimension $m$, that is $M = P/SO(m)$,
where $P$ is the bundle of oriented frames on the orbifold tangent bundle
(cf. section 1.4). Then the $C^*$ algebra of the orbifold $M$ is by
fiat the crossed product $C^*(M) = C(P) \rtimes SO(m)$, where $C(P)$
denotes the $C^*$ algebra of continuous functions on $P$. We will now
study some Morita equivalent descriptions of $C^*(M)$ that will be useful
for us later.
The following is one such, and is due to \cite{Far}.

\begin{prop} Let $M$ be a good orbifold, which is orbifold
covered by the smooth manifold $X$, i.e. $M = X/G$. Then
the $C^*$ algebras $C_0(X) \rtimes G$ and $C^*(M)$
are strongly Morita equivalent.
\end{prop}

In the two dimensional case, there is yet another $C^*$
algebra that is strongly Morita equivalent to the $C^*$
algebra of the orbifold. Let
$\Gamma$ be as before. Then $\Gamma$ acts freely
on $PSL(2, \R)$, and therefore the quotient space
$\Gamma\backslash PSL(2, \R) =P(g, \nu_1,\ldots,\nu_n)$ is a smooth compact
manifold,
with a right action of $SO(2)$ that is only infinitismally free.
The $C^*$ algebra of the hyperbolic orbifold
$\Sigma(g, \nu_1,\ldots,\nu_n)$ is by fiat the crossed product $C^*$
algebra
$$
C^*(\Sigma(g, \nu_1,\ldots,\nu_n)) =
C(P(g, \nu_1,\ldots,\nu_n)) \rtimes SO(2)
$$
cf. \cite{Co}. If $SO(2)$ did act freely on $P(g, \nu_1,\ldots,\nu_n)$
(which is the case when $\nu_1=\ldots =\nu_n=1$), then it is
known that
$C^*(\Sigma(g, \nu_1,\ldots,\nu_n)) $ and
$C(\Sigma(g, \nu_1,\ldots,\nu_n)) $ are strongly Morita equivalent
as $C^*$ algebras.

We shall next describe a natural algebra which is Morita equivalent
to the $C^*$ algebra of the orbifold $\Sigma(g, \nu_1,\ldots,\nu_n)$.
Now $\Gamma$ has a torsionfree subgroup $\Gamma_{g'}$ of finite index, such
that
the quotient $\Gamma_{g'}\backslash \J = \Gamma_{g'} \backslash PSL(2,\R)/SO(2)
= \Sigma_{g'}$ is a compact Riemann
surface of genus $g' =
1+ \frac{\#(G)}{2}(2(g -1) +(n - \nu))$
where $\#(G) = 1+ \sum_{j=1}^n (\nu_j-1)$ and  where
$\nu = \sum_{j=1}^n 1/\nu_j$, cf. Theorem 2.5, \cite{Sc} and the orbifold
Euler characteristic calculations in there. Then $G \to \Sigma_{g'}
\to \Sigma(g, \nu_1,\ldots,\nu_n)$ is a finite orbifold cover, i.e. a
ramified covering space, where $G = \Gamma_{g'}\backslash \Gamma$.

\begin{prop}
The $C^*$ algebras $C(\Sigma_{g'})\rtimes G$,  $C^*(\Sigma(g,
\nu_1,\ldots,\nu_n))$ and $C_0(\mathbb H)\rtimes \Gamma$
are strongly Morita equivalent to each other.
\end{prop}

\begin{proof} The strong Morita equivalence of the last two $C^*$
algebras is contained in the previous Proposition. Since strong Morita
equivalence is an equivalence relation, it suffices to prove that
the first two $C^*$ algebras are strongly Morita equivalent.
Let $\hat P = \Gamma_{g'} \backslash PSL(2,\R)$
where $SO(2)$ acts on $\hat P$ the right, and therefore commutes with the left
$G$ action on $\hat P$. Moreover, the actions of $G$ and $SO(2)$
on $\hat P$ are free, and therefore one can apply a theorem of Green,
\cite{Green}
which implies in particular that $C_o(G\backslash \hat P)\rtimes SO(2) $ and
$C_o(\hat P/SO(2))\rtimes G $
are strongly Morita equivalent, i.e.
$C_o(P(g, \nu_1,\ldots,\nu_n)\rtimes SO(2) $ and
$C_o(\Sigma_{g'}) \rtimes G $
are strongly Morita equivalent, proving the proposition.
\end{proof}

\subsection{Orbifold vector bundles and $K$-theory}
Because of the Morita equivalences of the last section 1.3, we can
give several alternate and equivalent descriptions of orbifold
vector bundles over orbifolds. Firstly, there is the description using
transition functions cf. \cite {Sat}, \cite{Kaw}. Equivalently,
one can view an orbifold vector bundle over an $m$ dimensional
orbifold $M$ as being an $SO(m)$ equivariant vector bundle over
the bundle $P$ of oriented frames of the orbifold tangent bundle.
In the case of a good orbifold $M$, which is orbifold covered
by a smooth manifold $X$. Let $G$ be the discrete group acting on $X$,
$G \to X \to M=X/G$. Then an orbifold vector bundle on $M$ is the
quotient ${\mathcal V}_M=G\backslash {\mathcal V}_X$ of a vector bundle over
$X$ by the $G$ action. Notice that an orbifold vector bundle is not a
vector bundle over $M$: in fact, the fibre at a singular point is
isomorphic to a quotient of a vector space by a finite group action.

The Grothendieck group of isomorphism classes of orbifold vector
bundles on the orbifold $M$ is called the
{\em orbifold $K$-theory} of $M$ and is denoted by $K^0_{orb}(M)$,
which by a result of \cite{BC}, \cite{Far} is canonically isomorphic
to $K_0(C^*(M))$. By the Morita
equivalence of section 1.3, one then has $K^0_{orb}(M) \cong K^0_{SO(m)}(P)$,
and by the Julg-Green theorem \cite{Ju}, \cite{Green}, the second group
is somorphic to $K_0(C(P) \rtimes SO(m))$. In the case when $M$ is
a good orbifold, by Proposition 1.1, one sees that $K^0_{orb}(M)
\cong K^0(C_0(X) \rtimes G) = K^0_G (X) $.

We will now be mainly interested in orbifold line bundles
over the hyperbolic 2-orbifolds. Let $G$ be the finite group
determined by the exact sequence
$1\to \Gamma_{g'}\to \Gamma\to G\to 1$. Then $G$ acts on $\Sigma_{g'}$ with
quotient the orbifold $\Sigma(g,
\nu_1,\ldots,\nu_n)$.

An orbifold line bundle ${\mathcal L}$ on $\Sigma(g,
\nu_1,\ldots,\nu_n)$ is given by
$$ {\mathcal L} =G\backslash (P\times_{SO(2)} \C),$$
where $P$ is a principal $SO(2)$-bundle on the smooth surface
$\Sigma_{g'}$. Notice that the $SO(2)$ and the $G$
actions commute, and are free on the total space $P$.
An orbifold line bundle has an associated Seifert fibred space $G
\backslash P$. A more explicit local geometric construction of ${\mathcal
L}$ is given in \cite{Sc}.
An orbifold line bundle ${\mathcal L}$ over a hyperbolic orbifold
$\Sigma(g,\nu_1,\ldots,\nu_n)$ is specified by the Chern class of the
pullback line bundle on the smooth surface $\Sigma_{g'}$,
together with the Seifert data. That is the pairs of numbers
$(\beta_j,\nu_j)$, where $\beta_j$ satisfies the following condition.
Given the exact sequence
$$ 1\to \Z \to \pi_1(P)\to \pi_1^{orb}(\Sigma(g,\nu_1,\ldots,\nu_n))\to 1,$$
let $\tilde C_j$ be an element of $\pi_1(P)$ that maps to the
generator $C_j$ of the orbifold fundamental group. Let $C$ be the
generator of the fundamental group of the fibre. Then we have
$C_j^{\nu_j}=1$ and $C^{\beta_j}=\tilde C_j^{\nu_j}$. The choice of
$\beta_j$ can be normalised so that $0<\beta_j <\nu_j$.

In more geometric terms, let $p_1, \ldots, p_n$ be the cone points of a
hyperbolic orbifold $\Sigma(g,\nu_1,\ldots,\nu_n)$. Let $\Sigma'$ be
the complement of the union of small disks around the cone points.
The orbifold line bundle induces a line bundle ${\mathcal L}'$
over the smooth surface with
boundary $\Sigma'$, trivialized over the
boundary components of $\Sigma'$. Moreover, the restriction of the
orbifold line bundle ${\mathcal L}$ over the
small disks $D_{p_i}$ around each cone point $p_i$ is obtained by
considering a surgery on the trivial product $\C\times D_{p_i}$
obtained by cutting open along a radius in $\C$ and gluing back after
performing a rotation on $D_{p_i}$ by an
angle $2\pi q/\nu_i$. With this notation the Seifert invariants are
$(q_i,\nu_i)$ with $\beta_i q_i\equiv 1$ $(mod \nu_i)$.

Thus, an orbifold line bundle has a finite set of singular fibres at
the cone points. The orbifold line bundle ${\mathcal L}$ pulls
back to a $G$-equivariant line bundle $\tilde {\mathcal L}$ over the
smooth surface $\Sigma_{g'}$  that orbifold covers
$\Sigma(g,\nu_1,\ldots,\nu_n)$. All the orbifold line bundles with
trivial orbifold Euler class, as defined in \cite{Sc}, lift to the
trivial line bundle on $\Sigma_{g'}$.

In \cite{Sc} the classification of Seifert-fibred
spaces is derived using the Seifert invariants, namely the Chern class
of the
line bundle $\tilde {\mathcal L}$, together with the Seifert data
$(\beta_j, \nu_j)$ of the singular fibres at the cone points $p_j$.
We show in the following that the Seifert invariants can be recovered
from the image of the Baum-Connes equivariant Chern character
\cite{BC}.

\subsection{Baum-Connes Chern character}

We have seen that the algebra $C^*(\Sigma(g,\nu_1,\ldots,\nu_n))$ is
strongly Morita equivalent to the cross product
$C(\Sigma_{g'})\rtimes G$. Therefore the relevant K-theory is
$$K_0(C(\Sigma_{g'})\rtimes G)=K^0_{SO(2)}(G\backslash \hat
P)=K^0_G(\Sigma_{g'}),$$
where $\hat P=\Gamma_{g'} \backslash PSL(2,\R)$.

We recall briefly the definition of delocalised equivariant cohomology
for a finite group action on a smooth manifold \cite{BC}. Let $G$ be a
finite group acting smoothly and properly on a compact smooth manifold $X$.
Let $M$ be the good orbifold $M=G\backslash X$. Given any $\gamma\in
G$, the subset
$X^\gamma$ of $X$ given by
$$ X^\gamma=\{ (x,\gamma)\in X\times G \mid \gamma x=x \} $$
is a smooth compact submanifold. Let $\hat X$ be the disjoint union of the
$X^\gamma$ for $\gamma\in G$. The complex $\Omega_G(\hat X)$ of
$G$-invariant de Rham
forms on $\hat X$ with coefficients in $\C$ computes the delocalised
equivariant cohomology $H^\bullet (X,G)$, which is $\Z_2$ graded by
forms of even and odd degree. The dual complex that computes
delocalised homology is
obtained by considering $G$-invariant de Rham currents on $\hat X$.
Thus we have
$$ H^\bullet (X,G)=H^\bullet (\Omega_G(\hat X),d)=H^\bullet (\hat X
/G,\C)$$
$$=H^\bullet (\hat X,\C)^G =\bigoplus_{\gamma\in G}
H^\bullet(X^\gamma,\C).$$
According to  \cite{BC}, Theorem 7.14, the delocalised equivariant
cohomology is isomorphic to the cyclic cohomology of the algebra
$C^\infty(X)\rtimes G$,
$$ H^0(X, G)\cong HC^{ev}(C^\infty(X)\rtimes G), $$
$$ H^1(X, G)\cong HC^{odd}(C^\infty(X)\rtimes G). $$

The Baum-Connes  equivariant Chern character
$$ ch_G : K^0_G(X)\to H^0(X,G) $$
is an isomorphism over the complex numbers. Equivalently, the
Baum-Connes  equivariant Chern character can be viewed as
$$ ch_G : K^0_{orb}(M)\to H^0_{orb}(M) $$
where the {\em orbifold cohomology} is by definition $ H^j_{orb}(M) =
H^j(X, G)$
for $j=0,1$.

In our case the delocalised equivariant cohomology and the Baum-Connes
Chern character have a simple expression.
Let $\Sigma_{g'}$ be the smooth surface that orbifold covers
$\Sigma(g,\nu_1,\ldots,\nu_n)$. Let $G$ be the finite group
$1\to \Gamma_{g'} \to \Gamma\to G\to 1$.
Let $G_{p_j}\cong \Z_{\nu_j}$ be
the stabilizer
of the cone point $p_j$ in $\Sigma(g,\nu_1,\ldots,\nu_n)$.
Then we have
$$ \Sigma_{g'}^\gamma=\left\{\begin{array}{lr}
\Sigma_{g'} & \quad\text{if} \quad\gamma=1; \\
\{p_j\} & \quad\text{if} \quad\gamma\in G_{p_j}\backslash \{ 1 \}; \\
\emptyset & \text{otherwise.} \end{array}\right. $$
Thus the delocalised equivariant cohomology and orbifold cohomology is given by
$$ H^0_{orb}(\Sigma(g,\nu_1,\ldots,\nu_n)) = H^0(\Sigma_{g'},
G)=H^0(\Sigma_{g'})\oplus H^2(\Sigma_{g'})\oplus
\C^{\sum_j (\nu_j-1)}, $$
where each $\C^{\nu_j-1}$ is given by $\nu_j-1$ copies of $H^0(p_j)$, and
$$ H^1_{orb}(\Sigma(g,\nu_1,\ldots,\nu_n)= H^1(\Sigma_{g'},
G)=H^1(\Sigma_{g'}). $$

Let ${\mathcal L}$ be an orbifold line bundle in
$K_0(C(\Sigma_g)\rtimes G)=K^0_G(\Sigma_g)$, and let $\tilde {\mathcal
L}$ be the corresponding line bundle over the surface $\Sigma_{g'}$. An
element $\gamma$ in the stabiliser $G_{p_j}$ acts on the restriction
of ${\mathcal L}|_{\Sigma^\gamma_g}={\mathcal L}|_{p_j}=\C$ as
multiplication by $\lambda(\gamma)=e^{2\pi i \beta_j/\nu_j}$.

Thus, the Baum-Connes Chern character of ${\mathcal L}$ is given by
$$ ch_G({\mathcal L})=(1,c_1(\tilde {\mathcal L}),
e^{2\pi i \beta_1/\nu_1}, \ldots, e^{2\pi i (\nu_1-1) \beta_1/\nu_1}, \ldots
e^{2\pi i \beta_n/\nu_n}, \ldots, e^{2\pi i (\nu_n-1) \beta_n}/\nu_n). $$

\begin{prop}
The Baum-Connes Chern character classifies orbifold line bundles over
$\Sigma(g,\nu_1,\ldots,\nu_n)$.
\end{prop}

\begin{proof}
According to \cite{Sc} the orbifold line bundles are classified by the
orbifold Euler number
$$ e(\Sigma(g,\nu_1,\ldots,\nu_n))=<c_1(\tilde {\mathcal
L}), [\Sigma_{g'}]>+ \sum_{j} \beta_j/\nu_j, $$
given in terms of the Chern number $<c_1(\tilde {\mathcal L}),
[\Sigma_{g'}]>$ and the Seifert invariants $(\beta_j,\nu_j)$.
\end{proof}

Notice that we have the isomorphism in $K$-theory,
$K^0_G(\Sigma_{g'})=K^0_{SO(2)}(G\backslash \hat P)$ and the Chern
character isomorphisms (with $\C$ coefficients)
$$ ch_G: K^0_G(\Sigma_{g'}) \to H^0(\Sigma_{g'}, G) \cong
HC^{ev}(C^\infty(\Sigma_{g'})\rtimes G) $$
and
$$ ch_{SO(2)} : K^0_{SO(2)}(\Gamma \backslash PSL(2,\R))\to
HC^{ev}(C^\infty(\Gamma \backslash PSL(2,\R))\rtimes SO(2)).$$
Moreover, we have an isomorphism
$$HC^\bullet (C^\infty(\Gamma \backslash PSL(2,\R))\rtimes
SO(2))\cong H^\bullet_{SO(2)} (\Gamma \backslash PSL(2,\R)). $$
Thus, we obtain
$$HC^{ev}(C^\infty(\Gamma \backslash PSL(2,\R))\rtimes SO(2)) \cong
HC^{ev}(C^\infty(\Sigma_{g'})\rtimes G)$$
with $\C$ coefficients, via the Chern character.

Thus orbifold line bundles on $\Sigma(g,\nu_1,\ldots,\nu_n)$ can be
described equivalently as $G$-equivariant line bundles over the
covering smooth surface $\Sigma_{g'}$, and again as
$SO(2)$-equivariant line bundles on $G\backslash \hat P$.

\begin{rems*}
With the notation used in the previous section, let $G$ be a
finite group acting smoothly and properly on a smooth compact
oriented manifold $X$.
There is a natural choice of a fundamental class $[X]_G \in H_0(X,G)$
in the delocalized equivariant homology of $X$, given by the
fundamental classes of each compact oriented smooth submanifold $X^\gamma$,
$[X]_G = \oplus_{\gamma\in G} [X^\gamma]$.
In the case of hyperbolic 2-orbifolds, the equivariant fundamental
class $[ \Sigma_{g'} ]_G$ is given by
$$ [ \Sigma_{g'} ]_G=[\Sigma_{g'}]\oplus_{j} [p_j]^{\nu_j-1}\in
H_2(\Sigma_{g'},\C)\oplus_j (H_0(p_j,\C))^{\nu_j-1}.$$
The corresponding equivariant Euler number $< ch_G({\mathcal L}),
[\Sigma_{g'}]_G >$ is obtained by evaluating
$$ < ch_G({\mathcal L}), [\Sigma_{g'}]_G >= <c_1(\tilde {\mathcal
L}), [\Sigma_{g'}]> + \sum_{j=1}^n \sum_{\gamma\in G_{p_j}\backslash
\{1 \}} \lambda(\gamma). $$
\end{rems*}

\subsection{Classifying space of the orbifold fundamental group}

Here we find it convenient to follow
Baum, Connes and Higson \cite{BC},
\cite{BCH}. Let $M$ be a
good orbifold, that is its orbifold universal cover $\widetilde M$
is a smooth manifold which has a
{\em proper} $\Gamma$-action, where $\Gamma$ denotes the orbifold
fundamental group of $M$. That is, the map
\begin{align*}
\widetilde M \times \Gamma & \to \widetilde M \times \widetilde M\\
(x, \gamma) & \to (x , \gamma x)
\end{align*}
is a proper map. The universal example for such a proper action
is denoted in \cite{BC}, \cite{BCH} by $\underline E\Gamma$. It is
universal in the
sense that there is a continuous $\Gamma$-map
$$
f: \widetilde M \to \underline E\Gamma
$$ which is unique upto $\Gamma$-homotopy, and moreover $\underline
E\Gamma$ itself is unique upto $\Gamma$-homotopy. The quotient
$\underline B\Gamma =  \Gamma\backslash \underline
E\Gamma$ is an orbifold. Just as $B\Gamma$ classifies isomorphism
classes of
$\Gamma$-covering spaces, it can be shown that $ \underline B\Gamma$
classifies  isomorphism classes of orbifold $\Gamma$-covering spaces.

\begin{examples*}
It turns out that if $\Gamma$ is a discrete subgroup of a connected
Lie group
$G$, then $\underline E\Gamma = G/K$, where $K$ is a maximal compact
subgroup. This is the main class of  examples that we are concerned
with in this paper.
\end{examples*}

Let $S\Gamma$ denote the set of all elements of $\Gamma$ which are of finite
order. Then $S\Gamma$ is not empty, since $1\in S\Gamma$. $\Gamma$
acts on $S\Gamma$ by conjugation, and let $F\Gamma$ denote
the associated permutation module over $\C$, i.e.
$$
F\Gamma = \left\{\sum_{\alpha \in \S\Gamma} \lambda_\alpha [\alpha]
\, \Big|\,
\lambda_\alpha \in \C \quad \text{and}\quad \lambda_\alpha = 0
\quad \text{except for
a finite number of}\quad\alpha\right\}
$$

Let $C^k(\Gamma: F\Gamma) $ denote the space of all antisymmetric
$F\Gamma$-valued $\Gamma$-maps on $\Gamma^{k+1}$, where $\Gamma$ acts
on  $\Gamma^{k+1}$ via the diagonal action. The coboundary map is
$$
\partial c(g_0, \ldots, g_{k+1}) =
\sum_{i=0}^{k+1}
(-1)^i c(g_0, \ldots,\hat g_i \ldots g_{k+1})
$$
for all $c \in C^k(\Gamma: F\Gamma)$ and
where $\hat g_i $ means that $g_i$ is omitted.
The cohomology
of this complex is the group cohomology of $\Gamma$
with coefficients in $F\Gamma$,
$H^k(\Gamma: F\Gamma) $, cf. \cite{BCH}. They also show that
$H^k(\Gamma: F\Gamma) \cong H^j(\Gamma, \C)\oplus_m H^k(Z(C_m):\C)$, where
$S\Gamma = \{1, C_m | m = 1, \ldots \}$ and the isomorphism is canonical.

Also, for any Borel measurable $\Gamma$-map
$\mu: \underline E\Gamma \to \Gamma$, there is an induced map
on cochains
$$
\mu^*: C^k(\Gamma: F\Gamma) \to C^k(\underline E\Gamma: \Gamma)
$$
which induces an isomorphism on cohomology,
$\quad \mu^* : H^k(\Gamma: F\Gamma)
\cong H^k(\underline E\Gamma: \Gamma)\quad$
\cite{BCH}. Here $H^j(\underline E\Gamma, \Gamma)$
denotes the $\Z$-graded (delocalised) equivariant cohomology
of $\underline E\Gamma$, which is a refinement of what was
discussed earlier, and which is
defined in \cite{BCH} using sheaves (and cosheaves), but we  will not
recall the definition here.

Let $M$ be a good orbifold with orbifold fundamental group $\Gamma$.
We have seen that the universal orbifold cover $\widetilde M$
is classified by a
continuous map $f: M \to \underline B\Gamma$, or equivalently
by a $\Gamma$-map $f: \widetilde M \to \underline E\Gamma$.
The induced map is $f^* :
H^j_{orb}(\underline B\Gamma, \C)  \equiv
H^k(\underline E\Gamma: \Gamma)
\to  H^k(\widetilde M: \Gamma) \equiv H^k_{orb}( M: \C)$ and therefore
in particular one has $f^*([c]) \in  H^k_{orb}( M: \C)$
for all $[c] \in H^k(\Gamma: \C)$.
This can be expressed on the level of cochains by a easily
modifying
the procedure in \cite{CM}, and we refer to
\cite{CM} for further details.

\section{ Twisted higher index theorem}

In this section, we will define the higher twisted index of an elliptic
operator on a good orbifold, and establish a cohomological formula for
any cyclic trace arising from a group cocycle, and
which is applied to the twisted higher index. We adapt the strategy
and proof in \cite{CM} to our context.

\subsection{Construction of the parametrix and and the index map}
For basic material on orbifolds, see \cite{Sc} and references there.
Let $M$ be a compact, good orbifold, that is, the universal cover
$\Gamma \to \M \to M$ is a smooth manifold
and we will assume as before that there is a $(\Gamma, \bar\sigma)$-action
on $L^2(\M)$ given by $T_\gamma = U_\gamma \circ S_\gamma \, \forall \gamma
\in \Gamma$. Let $\E, \ \F$ be Hermitian vector bundles on $M$ and let
$\E, \ \F$ be the corresponding lifts to $\Gamma$-invariants Hermitian vector
bundles on $\M$. Then there are induced $(\Gamma, \sigma)$-actions
on $L^2(\M, \E)$ and $L^2(\M, \F)$ which are also given by
$T_\gamma = U_\gamma \circ S_\gamma \, \forall \gamma \in \Gamma$.

Now let $D : L^2(\M, \E) \to L^2(\M, \F)$ be a 1st order
$(\Gamma, \bar\sigma)$-invariant elliptic operator. Let $U\subset \M$
be an open subset that contains the closure of a fundamental domain
for the $\Gamma$-action on $\M$. Let $\psi \in C^\infty_c(\M)$ be a
compactly supported smooth function such that $\text{supp}(\psi) \subset
U$, and
$$
\sum_{\gamma \in \Gamma} \gamma^*\psi = 1.
$$
Let $\phi \in C^\infty_c(\M)$ be a
compactly supported smooth function such that $\phi=1$
on $\text{supp}(\psi)$.

Since $D$ is elliptic, we can construct a parametrix $J$ for it on the open
set $U$ by standard methods,
$$
JDu = u - Hu \qquad \forall u \in C^\infty_c(U, \E|_U)
$$
where $H$ has a smooth Schwartz kernel.
Define the pseudodifferential operator $Q$ as
$$
Q = \sum_{\gamma \in \Gamma} T_\gamma \phi J \psi T_\gamma^*\eqno{(1)}
$$
We compute,
$$
QD w = \sum_{\gamma \in \Gamma} T_\gamma \phi J \psi D T_\gamma^* w \qquad
\forall w \in C^\infty_c(\M, \E), \eqno{(2)}
$$
since $T_\gamma D = D T_\gamma \quad \forall\gamma\in \Gamma$. Since $D$ is
a 1st
order operator, one has
$$
D(\psi w ) = \psi Dw + (D\psi) w
$$
so that $(2)$ becomes
$$
=\sum_{\gamma \in \Gamma} T_\gamma \phi J D \psi  T_\gamma^* w
- \sum_{\gamma \in \Gamma} T_\gamma \phi J (D\psi)  T_\gamma^* w.
$$
Using $(1)$, the expression above becomes
$$
= \sum_{\gamma \in \Gamma} T_\gamma \psi  T_\gamma^* w
- \sum_{\gamma \in \Gamma} T_\gamma \phi H \psi  T_\gamma^* w
- \sum_{\gamma \in \Gamma} T_\gamma \phi J (D \psi)  T_\gamma^* w.
$$
Therefore $(2)$ becomes
$$
QD = I - R_0
$$
where
$$
R_0 = \sum_{\gamma \in \Gamma} T_\gamma \left(\phi H \psi
+  J (D \psi) \right) T_\gamma^*
$$
has a smooth Schwartz kernel.
It is clear from the definition that one has $T_\gamma Q = Q T_\gamma$
and  $T_\gamma R_0 = R_0 T_\gamma \quad \forall\gamma\in \Gamma$.
Define
$$
R_1 = ^tR_0 + D R_0 ^tQ - DQ(^t R_0).
$$
Then  $T_\gamma R_1 = R_1 T_\gamma \quad \forall\gamma\in \Gamma$,
$R_1$ has a smooth Schwartz kernel and satisfies
$$
DQ = I - R_1.
$$
Summarizing, we have the following

\begin{prop}
Let $M$ be a compact, good orbifold and $\Gamma \to \M \to M$
be the universal orbifold
covering space. Let $\mathcal E, \ \mathcal F$ be Hermitian vector bundles on
$M$ and let
$\E, \ \F$ be the corresponding lifts to $\Gamma$-invariants Hermitian vector
bundles on $\M$. We will assume as before that there is a $(\Gamma,
\bar\sigma)$-action
on $L^2(\M)$ given by $T_\gamma = U_\gamma \circ S_\gamma \, \forall \gamma
\in \Gamma$. Then there are $(\Gamma, \bar\sigma)$-actions
on $L^2(\M, \E)$ and $L^2(\M, \F)$ which are also given by
$T_\gamma = U_\gamma \circ S_\gamma \, \forall \gamma \in \Gamma$.

Now let $D : L^2(\M, \E) \to L^2(\M, \F)$ be a 1st order
$(\Gamma, \sigma)$-invariant elliptic operator. Then there
is an almost local $(\Gamma, \bar\sigma)$-invariant
elliptic pseudodifferential operator $Q$
and $(\Gamma, \bar\sigma)$-invariant smoothing operators
$R_0, \ R_1$ which satisfy
$$
QD = I - R_0 \qquad {\text{and}} \qquad DQ = I - R_1.
$$
\end{prop}

Define
$$
   e(D) = \begin{pmatrix} {R}_0^2 & ({R}_0+
    {R}^2_0) {Q} \\ {R}_1{D} &
      1-{R}^2_1 \end{pmatrix}
$$
then $e(D)\in M_2\bigl(C^\infty_c(\mathcal{G},\sigma)\bigr)$ is an
{\em idempotent}, where  $C^\infty_c(\mathcal{G},\sigma)$ is as
defined in section 1. The
$C^\infty_c(\mathcal{G},\sigma)$-{\em index map} is by fiat
\[
   \Ind_\sigma(D)=[e(D)] - [e_0]\in
K_0(C^\infty_c\bigl(\mathcal{G},\sigma)\bigr)
\]
where $e_0$ is the idempotent
$$
e_0 = \begin{pmatrix} 0 & 0 \\ 0 &
      1 \end{pmatrix}.
$$
It is not difficult to see that $ \Ind_\sigma(D)$ is independent of the choice
of $(\Gamma, \sigma)$-invariant parametrix $Q$ that is needed in its
definition.

\begin{lemma}
Consider the right action of the algebra $C^\infty(M) \otimes \C(\Gamma,
\sigma)$ on the
vector space of compactly supported smooth functions on $\M$,
$C^\infty_c(\M)$ given by
$$
(\xi(f\otimes T_\gamma))(x) = f(p(x)) (T_\gamma\xi)(x) \qquad
\forall \xi \in
C^\infty_c(\widetilde M), \forall x\in \widetilde M \,
\text{and }\, \forall \gamma \in \Gamma,
$$
where $(T_\gamma\xi)(x) = (U_\gamma S_\gamma\xi)(x)$, \
$(U_\gamma\xi)(x) =\xi(\gamma x)$ and
$(S_\gamma\xi)(x) = e^{i\phi_\gamma(x)} \xi(x)$ for all
$\gamma\in \Gamma$. Then
$C^\infty_c(\M)$ is a finite projective module over
$C^\infty(M) \otimes \C(\Gamma, \sigma)$.
\end{lemma}

\begin{proof}
Let $B = C^\infty(M) \otimes \C(\Gamma, \sigma)$ and $E =
C^\infty_c(\M)$. Let $\{V_i\}_{i=1}^N$ be a finite open
cover of $M$, $\beta_i : {V_i}
\to \widetilde M$ be a smooth section of the orbifold
covering
$p: \widetilde M \to M$  for $j=1,\ldots N$. Let $\{\chi_i
\}_{i=1}^N$ be a partition of unity of $M$ that is subordinate to the open
cover $\{V_i\}_{i=1}^N$ such that the functions $\{\chi_i^{1/2}\}_{i=1}^N$
are smooth
on $M$. For $\xi \in E$ , define
$$
U \in \text{Hom}_B(E, B^N)
$$
by
$$
(U\xi)(x, \gamma) = (\chi_1^{1/2}(x) \xi(\gamma^{-1}\beta_1 (x)),
 \ldots, \chi_N^{1/2}(x) \xi(\gamma^{-1}\beta_N (x)))
\qquad \forall x\in M \forall \gamma\in \Gamma.
$$
For $\eta \in B^N$ and $\eta = (\eta_1,\ldots, \eta_N)$, define
$$
U^* \in \text{Hom}_B( B^N, E)
$$
by
$$
(U^*\eta)(x) = \sum_{j=1}^N \chi_j^{1/2}(p(x))
\eta_j(x, \gamma(\beta_j(p(x)), x)) \qquad \forall x\in \widetilde M.
$$
where $\gamma= \gamma(\beta_j(x), \tilde x)$ is the unique element in $\Gamma$
satisfiying $\gamma x = \beta_j(p(x))$.
It is straightforward to verify that
$U^*U = 1_E \in \text{Hom}_B( E, E) $ and therefore one sees that
$p_\sigma = UU^* \in \text{Hom}_B( B^N,B^N) = M_N(B)$ is an idempotent.
Explicitly, one has
$$
(p_\sigma)_{jk} = \chi_j^{1/2}\chi_k^{1/2} \otimes
\beta_j\beta_k^{-1}\qquad \text{for}\, 1\le j,k \le N.
$$
Therefore $U: E \to p_\sigma B^N$ is an isomorphism of $B$-modules, and
therefore
of $\C(\Gamma, \sigma)$-modules, where $\C(\Gamma, \sigma)$ is identified with
the subalgebra $1\otimes \C(\Gamma, \sigma)$ of $B = C^\infty(M)
\otimes \C(\Gamma, \sigma)$.
\end{proof}

Using this lemma, we can define a homomorphism $\theta_\sigma$, defined by
$\theta_\sigma(T) = U T U^*$, which transforms $(\Gamma, \sigma)$-invariant
linear operators on  $C^\infty_c(\M)$ into $\C(\Gamma, \sigma)$-linear
endomorphisms
of $p_\sigma B^N$. Since $B = C^\infty(M) \otimes \C(\Gamma, \sigma)$, an
arbitrary
$\C(\Gamma, \sigma)$-linear endomorphism of $p_\sigma B^N$ is given by a matrix
$\left(S_{i,j}\right)_{1\le i, j\le N}$, where $ S_{i,j} \in
\End(C^\infty(M)) \otimes
\C(\Gamma, \sigma)$.

Let $\mathcal{R}_M$ denote the algebra of all smoothing operators on
$C^\infty(M)$.

\begin{prop}
The homomorphism $\theta_\sigma$ above, maps the algebra
$C^\infty_c(\mathcal{G},\sigma)$
into $M_N(\mathcal{R}_M \otimes \C(\Gamma, \sigma))$.
\end{prop}

\begin{proof} The map is explicitly given by the equality
$$
\theta_\sigma(k)_{ij}(x,y,\gamma) = \chi_i(x)^{1/2}\chi_j(y)^{1/2}
k(\beta_j(x), \gamma\beta_j(y))$$
$\forall k \in C^\infty_c (\M\times_\Gamma \M) =
C^\infty_c(\mathcal{G},\sigma),
\ \forall x,y \in M, \ \forall \gamma\in \Gamma$,
and therefore $\theta_\sigma(k)_{ij}(\gamma)$ is clearly a smoothing operator.
\end{proof}

\begin{prop}
The homomorphism $\theta_\sigma$ induces a homomorphism
$$
K_0(C^\infty_c(\mathcal{G},\sigma)) \to K_0(\mathcal{R}_M \otimes
\C(\Gamma, \sigma))
$$
that is independent of the choice of $(U_j, \beta_j, \chi_j)_{1\le j\le N}$.
\end{prop}

\begin{proof}
Suppose that $U' \in \text{Hom}_B(E, B^N)$ is another choice of $U$. At the
expense
of replacing $N$ by $2N$, we may assume that there is $W \in GL_N(B)$ such
that
$$
U' = W U.
$$
Then ${\theta'}_\sigma(T) = U' T {U'}^* = W \theta_\sigma(T) W^{-1} \quad
\forall
T \in C^\infty_c(\mathcal{G},\sigma)$, where $W$ is viewed as a multiplier
of the
algebra $M_N(\mathcal{R}_M \otimes \C(\Gamma, \sigma))$. Therefore both
$\theta'_\sigma$ and $\theta_\sigma$ induce the same map on $K$-theory.
\end{proof}

The following lemma is an immediate generalisation of a result in \cite{Co}

\begin{lemma}
Let $M$ be a compact orbifold of positive dimension. Then there is a canonical
isomorphism
$$
\rho : \mathcal{R}_M \to \mathcal{R}
$$
which is unique upto inner automorphisms of $\mathcal{R}$.
\end{lemma}

Combining the Proposition 2.4 and Lemma 2.5 above, we obtain a canonical
homomorphism
in $K$-theory
$$
J: K_0(C^\infty_c(\mathcal{G},\sigma)) \to K_0(\mathcal{R}(\Gamma, \sigma))
$$
where $\mathcal{R}(\Gamma, \sigma) = \mathcal{R} \otimes \C(\Gamma, \sigma)$,
and $J = (\rho\otimes 1) \otimes \theta_\sigma$.

\begin{defn*}
The $(\Gamma, \sigma)$-{\em index} of a $(\Gamma, \sigma)$-invariant
elliptic operator
$D : L^2(\M, \E) \to L^2(\M, \F)$ is defined as
$$
\Ind_{(\Gamma, \sigma)}(D) = J(\Ind_\sigma(D)) \in K_0(\mathcal{R}(\Gamma,
\sigma))
$$
\end{defn*}

\subsection{Heat kernels and the index map} Let $D$ be as before, and for
$t>0$,
using the standard off-diagonal estimates for the heat kernel
recall that the heat kernels $e^{-tD^*D}$ and $e^{-tDD^*}$ are elements in the
$M_N\bigl(C^\infty_c(\mathcal{G},\sigma)\bigr)$ for some $N$ large enough.
Define the idempotent $e_t(D)\in
M_{2N}\bigl(C^\infty_c(\mathcal{G},\sigma)\bigr)$
as follows
$$
   e_t(D) = \begin{pmatrix} e^{-t D^*D} & e^{-t/2 D^*D}\frac{(1-e^{-t D^*D})}
   {D^*D} D \\
e^{-t/2 D D^*}{D} &
      1- e^{-t DD^*} \end{pmatrix}
$$
where $f$ is a smooth, even function on $\R$ which satisfies $f(x) x^2 =
e^{-x^2} (1-e^{-x^2} )$.

The relationship with the idempotent $e(D)$ constructed earlier will be
explained now. Define for $t>0$,
$$
Q_t = \frac{\left( 1-e^{-t/2 D^*D}\right)}{D^*D} D^*
$$
Then one easily verifies that $Q_t D = 1-e^{-t/2 D^*D} = 1- R_0(t)$
and $D Q_t  = 1-e^{-t/2 D D^*} = 1- R_1(t)$. That is, $Q_t$ is a parametrix
for $D$ for all $t>0$. Therefore one can write
$$
e_t(D) = \begin{pmatrix} R_0(t)^2 & (R_0(t) + R_0(t)^2) Q_t\\
R_1(t){D} &
     1- R_1(t)^2 \end{pmatrix}.
$$
In particular, one has for $t>0$
$$
\Ind_\sigma(D) = [e_t(D)] - [e_0] \in K_0(C^\infty_c(\mathcal{G},\sigma)).
$$

\subsection{Twisting an elliptic operator}
We will discuss elliptic operators only on good orbifolds, and refer to
\cite{Kaw} for the general case.
Let $M$ be a good orbifold, that is the universal orbifold cover
$\widetilde M$ of $M$ is a smooth manifold.
Let $\W\to \M$ be a $\Gamma$-invariant Hermitian vector bundle
over $\widetilde M$.
Let $D$ be a 1st order elliptic differential operator on $M$,
$$
D: L^2 (M, \mathcal E) \to L^2(M, \mathcal F)
$$
acting on $L^2$ orbifold sections of the orbifold vector bundles
$\mathcal E, \mathcal F$ over $M$. By fiat, $D$ is a $\Gamma$-equivariant
1st order
elliptic differential operator $\widetilde D$ on the smooth manifold
$\widetilde M$,
$$
\widetilde D :L^2 (\widetilde M, \E) \to L^2(\widetilde M, \F).
$$
Given any connection
$\nabla^{\W}$ on $\W$ which is compatible with the $\Gamma$ action and
the Hermitian metric, we wish to define an extension of the
elliptic operator $\widetilde D$, to act
on sections of $\E\otimes \W$, $\F\otimes \W$.
\[
\widetilde D\otimes\nabla^{\W}:\Gamma(M,\E\otimes \W)\to
\Gamma(M,\F\otimes \W)
\]
and we want it to satisfy the following property:\ If $\W$
is a trivial bundle,
and $\nabla^0$ is the trivial connection on $\W$, then for
$u \in \Gamma(\M,\E),\
h \in \Gamma(\M, \W)$ such that $\nabla^0 h = 0$,
\[
(\widetilde D\otimes\nabla^0)(u\otimes h)=
(\widetilde D u)\otimes h
\]
To do this, define a morphism
\begin{align*}
&\quad S=S_D\,:\,\E\otimes T^* \M\to\F\\ &\quad S(u\otimes df) =
\widetilde D(f u)-f{\widetilde D}u \end{align*}
for $f\in C^\infty(\M)$ and $u\in\Gamma(M,\E)$. Then $S$ is a tensorial.
Consider $S=S\otimes 1:\E\otimes T^* \M \otimes \W\to\F\otimes \W$
defined by \[
S(u\otimes df\otimes e)=S(u\otimes df)\otimes e \]
for $u, f$ as before and $e\in\Gamma(M,\W)$.\\
Recall that a connection $\nabla^{\W}$ on $\W$ is a derivation
\[
\nabla^{\W}:\Gamma(\M,\W)\to\Gamma(\M,T^* \M\otimes \W) \]
Define $\widetilde D\otimes\nabla^{\W}$ as
\[
(\widetilde D\otimes\nabla^{\W})(u\otimes e)=(D
u)\otimes e+S(u\otimes\nabla^{\W} e) \]
Then $\widetilde D\otimes\nabla^{\W}$ is a $1^{\text{ st }}$
order elliptic operator.

\subsection{Group cocycles and cyclic cocycles}

Using the pairing theory of cyclic cohomology and $K$-theory,
due to \cite{Co}, we will pair the $(\Gamma, \sigma)$-index of
a $(\Gamma, \sigma)$-invariant elliptic operator $D$ on $\widetilde M$
with certain cyclic cocycles on ${\mathcal R}(\Gamma, \sigma)$. The
cyclic cocycles that we consider come from normalised group
cocycles on $\Gamma$. More precisely,
given a normalised group cocycle $c\in Z^k(\Gamma, \mathbb C)$, $
k=0,\ldots, \dim M$,
we define a cyclic cocycle
$\tr_c$ of dimension $k$ on the twisted group ring
$\mathbb C(\Gamma, \sigma)$, which is
given by
$$
\tr_c(a_0\delta_{g_0}, \ldots,a_k\delta_{g_k}) = \left\{\begin{array}{l}
a_0\ldots a_k c(g_1,\ldots,g_k)
\tr( \delta_{g_0}\delta_{g_1} \ldots\delta_{g_k}) \quad\text{if} \,\, g_0
\ldots g_k =1\\
\\
0 \qquad \text{otherwise.}
\end{array}\right.
$$
where $a_j \in \C$ for $j=0,1,\ldots k$.
To see that this is a cyclic cocycle on $\C(\Gamma, \sigma)$, we first define
as done in \cite{Ji}, the twisted differential graded algebra
$\Omega^\bullet(\Gamma, \sigma)$
as the differential graded algebra of finite linear combinations of
symbols
$$
g_0 dg_1\ldots dg_n \qquad g_i\in \Gamma
$$
with module structure and differential given by
\begin{align*}
(g_0 dg_1\ldots dg_n)g
&= \sum_{j=1}^n(-1)^{n-1} \sigma(g_j, g_{j+1}) g_0dg_1\ldots d(g_jg_{j+1})
\ldots dg_n  dg\\
& + (-1)^n \sigma(g_n, g) g_0 dg_1\ldots d(g_ng)\\
d(g_0 dg_1\ldots dg_n) & = dg_0 dg_1\ldots dg_n
\end{align*}
We now recall normalised group cocycles. A group $k$-cocycle is
a map
$h:\Gamma^{k+1} \to \C$ satisfying the identities
\begin{align*}
h(g g_0, \ldots g g_k) & = h(g_0, \ldots g_k)\\
0 & = \sum_{i=0}^{k+1} (-1)^i h(g_0, \ldots, g_{i-1},g_{i+1}\ldots, g_{k+1})
\end{align*}
Then a normalised group $k$-cocycle $c$ that is
associated to such an $h$ is given by
$$
c(g_1, \ldots, g_k) = h (1, g_1, g_1 g_2, \ldots, g_1\ldots g_k)
$$
and it is defined to be zero if either $g_i = 1$ or if $g_1\ldots g_k = 1$.
Any normalised group cocycle $c \in Z^k(\Gamma, \C)$
determines a $k$-dimensional
cycle via the following closed graded trace on
$\Omega^\bullet(\Gamma, \sigma)$
$$
\int g_0 dg_1\ldots dg_n =
\left\{\begin{array}{l}
 c(g_1,\ldots,g_k)
\tr( \delta_{g_0}\delta_{g_1} \ldots\delta_{g_k}) \quad\text{if} \, \, n=k
\quad\text{and}\,\, g_0
\ldots g_k =1\\
\\
0 \qquad \text{otherwise.}
\end{array}\right.
$$

The higher cyclic trace $\tr_c$ is by fiat this closed graded
trace.

\subsection{Twisted higher index theorem-
the cyclic cohomology version}

Let $M$ be a compact orbifold of dimension $n=4\ell$.
Let $\Gamma\to\widetilde{M}\overset{p}{\to}M$ be the universal
cover of $M$ and the orbifold fundamental group is $\Gamma$.
Let $D$ be an elliptic 1st order
operator on $M$ and $\widetilde{D}$ be the lift of
$D$ to $\widetilde{M}$,
\[
\widetilde{D}\,:\,L^2(\widetilde{M},
\E)\to
L^2(\widetilde{M},\F). \]
Note that $\widetilde{D}$ commutes with the $\Gamma$-action on $\widetilde{M}$.

Now let $\omega$ be a closed 2-form on $M$ such that $\widetilde{\omega}
=p^* \omega=d\eta$ is \emph{exact}. Define $\nabla=d+\,i\eta$.
Then $\nabla$ is a Hermitian connection on the trivial line bundle
over $\widetilde{M}$, and the curvature of $\nabla,\ (\nabla)^2=i\,
\widetilde{\omega}$. (Here $s\in\mathbb{R}$.) Then $\nabla$ defines
a projective action of $\Gamma$ on $L^2$ spinors as follows:\\
Firstly, observe that since $\widetilde{\omega}$ is $\Gamma$-invariant,
$0=\gamma^*\widetilde{\omega}-\widetilde{\omega}=d(\gamma^*\eta-\eta)\
\forall\gamma \in\Gamma$. So $\gamma^*\eta-\eta$ is a closed 1-form on the
simply connected manifold $\widetilde{M}$, therefore \[
\gamma^*\eta-\eta=d\phi_\gamma\quad\forall\gamma\in\Gamma \]
where $\phi_\gamma$ is a smooth function on $\widetilde{M}$ satisfying in
addition,
\begin{itemize}
\item $\phi_\gamma(x)+\phi_{\gamma'}(\gamma x)-\phi_{\gamma'\gamma}(x)$
is independent of $x\in\widetilde{M}\ \forall \gamma,\gamma'\in\Gamma$;
\item $\phi_\gamma(x_0)=0$ for some $x_0\in\widetilde{M}\ \forall\in\Gamma$.
\end{itemize}
Then $\sigma(\gamma,\gamma')=\exp(si\phi_\gamma(\gamma'\cdot x_0))$
defines a multiplier on $\Gamma$ i.e. $\sigma:\Gamma\times\Gamma\to U(1)$
satisfies the following identity for all $\gamma, \gamma', \gamma'' \in \Gamma$
\[ \sigma(\gamma,\gamma')\sigma(\gamma,\gamma'\gamma'')=
\sigma(\gamma\gamma',\gamma'') \sigma(\gamma',\gamma'') \]
For $u \in L^2(\widetilde{M},
\E)$, let $
S_\gamma u  =  e^{i\phi_\gamma}u$ and $U_\gamma u =
\gamma^*u$
and $T_\gamma=U_\gamma {\text{\small o}} S_\gamma$ be the composition.
Then $T$ defines a projective $(\Gamma,\sigma)$-action on $L^2$-spinors,
i.e. \[
T_\gamma T_{\gamma'}=\sigma(\gamma,\gamma')T_{\gamma\gamma'}. \]
Consider the {\em twisted elliptic operator} on $\widetilde{M}$, \[
\widetilde{D}\otimes\nabla\,:\,L^2(\widetilde{M},
\E)\to
L^2(\widetilde{M},\F) \]
Then $\widetilde{D}\otimes\nabla$ no longer commutes with
$\Gamma$, but it does commute with the projective $(\Gamma,\sigma)$ action.
Let $P_+, P_-$ be the orthogonal projections onto the nullspace of
$\widetilde{D}\otimes\nabla$ and
$(\widetilde{D}\otimes\nabla)^*$ respectively since \[
(\widetilde{D}\otimes\nabla)\ P_+ = 0 \qquad {\text{and}}
\qquad (\widetilde{D}\otimes\nabla)^*\ P_- = 0\]
By elliptic regularily, it follows that the Schwartz (or integral) kernels of
$P_\pm$ are smooth. Since $\widetilde{D}\otimes\nabla$ and its adjoint
commutes with the $(\Gamma,\sigma)$ action, one has $$
e^{i\phi_\gamma(x)} P_\pm(\gamma x,\gamma y)\,e^{-i\phi_\gamma(y)}
= P_\pm(x,y)\quad\forall\gamma\in\Gamma.
$$
In particular, $P_\pm(x,x)$ is a $\Gamma$-invariant function on
$\widetilde{M}$.
One can define the {\em von Neumann trace} as Atiyah did in the untwisted case
\[
\tr\left(P_\pm\right) =
\int_M\,\tr \left(P_\pm(x,x)\right)\,dx. \]
The $L^2$-index is by definition
$$
\Index_{L^2} (\widetilde{D}\otimes\nabla)
= \tr(P_+)-\tr(P_-).
$$

To describe the next theorem, we will briefly review some material on
characteristic classes for orbifold vector bundles. Let $M$ be a good
orbifold, that is the universal orbifold cover $\Gamma \to \widetilde M
\to M$ of $M$
is a smooth manifold. Then the orbifold tangent bundle $TM$ of $M$,
can be
viewed as the $\Gamma$-equivariant bundle $T\widetilde M$ on $\widetilde M$.
similar comments apply to the orbifold cotangent bundle $T^*M$ and more
generally,
any orbifold vector bundle on $M$. It is then clear that choosing
$\Gamma$-invariant connections on the $\Gamma$-invariant vector bundles
on $\widetilde M$, one can define the Chern-Weil representatives of
the characteristic classes of the $\Gamma$-invariant vector bundles
on $\widetilde M$. These characteristic classes are $\Gamma$-invariant
and so define cohomology classes on $M$. For further details, see \cite{Kaw}.

\begin{thm}
Let $M$ be a compact, even dimensional, good
orbifold, $\Gamma$
be its orbifold fundamental group, $\widetilde D$ be a
$\Gamma$-invariant elliptic differential operator on $\widetilde M$, where
$\Gamma \to \widetilde M\to M$ is the universal orbifold cover of $M$.
Then for any group cocycle $c\in Z^{2q}(\Gamma)$, $q=0,2$, one has
$$
\Ind_{c, \Gamma, \sigma}(\widetilde D\otimes \nabla)
= \frac{q !} {(2\pi i)^q (2q!)} \left<Td(M)\cup ch(symb(D))
\cup f^*(\phi_c) \cup e^{\omega}, [T^*M]   \right>
$$
where $Td(M)$ denotes the Todd characteristic class of
the complexified orbifold tangent bundle of $M$
which is pulled back to the orbifold cotangent bundle $T^*M$,
$ch(symb(D))$ is the Chern character of the symbol of the operator
$D$, $\phi_c$ is the Alexander-Spanier cocycle on $\underline B\Gamma$
that corresponds to the group cocycle $c$ and $f: M
\to \underline B\Gamma$ is the map that classifies the
orbifold universal cover  $ \M \to M$, cf. sections 1.7 and 2.5.

\end{thm}

\begin{proof}
Choose a bounded, almost everywhere smooth Borel cross-section
$\beta : M\to \M$, which can then be used to define
the Alexander-Spanier cocycle $\phi_c$ corresponding to
$c\in Z^{2q}(\Gamma)$,  and such that $[\phi_c] = [c]
\in H^{2q}(\Gamma)$. As in Propositon 2.1, there is a
is an almost local $(\Gamma,
\sigma)$-invariant parametrix $Q$ of $\widetilde D\otimes \nabla$
and $(\Gamma, \sigma)$-invariant smoothing operators
$R_0, \ R_1$ which satisfy
$$
Q(\widetilde D\otimes \nabla) = I - R_0 \qquad {\text{and}}
\qquad (\widetilde D\otimes \nabla) Q = I - R_1.
$$
Then as before, one can then construct the index idempotent
$$
e(\widetilde D\otimes \nabla)
= \begin{pmatrix} {R}_0^2 & ({R}_0+
    {R}^2_0) {Q} \\ {R}_1{(\widetilde D\otimes \nabla)} &
     1- {R}^2_1 \end{pmatrix}
\in M_2\bigl(C^\infty_c(\mathcal{G},\sigma)\bigr)
$$
and the $C^\infty_c(\mathcal{G},\sigma)$-{\em index map}
is by fiat
\[
   \Ind_\sigma(\widetilde D\otimes \nabla)=
[e(\widetilde D\otimes \nabla)] - [e_0]\in
K_0(C^\infty_c\bigl(\mathcal{G},\sigma)\bigr).
\]
where $e_0$ is the idempotent
$$
e_0 = \begin{pmatrix} 0 & 0 \\ 0 &
      1 \end{pmatrix}.
$$
Let $R_t =e_t(\widetilde D\otimes \nabla) - e_0$ and $\theta_\sigma:
C^\infty_c(\mathcal{G},\sigma) \to \C(\Gamma,\sigma)\otimes {\mathcal L}^2
$ be the homomorphism obtained from the section $\beta$. Then one
has for $t>0$
$$
\Ind_{1, \Gamma, \sigma}(\widetilde D\otimes \nabla)
= \tr (\theta_\sigma(R_t))
$$
if $c = 1$ and proceeds exactly as in the case of
Atiyah's $L^2$ index theorem
for covering spaces. There is a standard reduction to the case when
$D = {{\not\partial}^\pm \otimes\nabla^{\mathcal E}} =
{\not\partial}^\pm_{\mathcal E} $.
Let $k^\pm(t,x,y)$ denote the heat kernel of the lifted Dirac operators
$(\widetilde{\not\partial}^\pm_{\mathcal E} \otimes \nabla)^2$
on the universal cover of $M$, and $P^\pm(x,y)$
the smooth Schwartz kernels of the orthogonal projections $P^\pm$ onto the
nullspace
of $\widetilde{\not\partial}^\pm_{\mathcal E}\otimes \nabla^s$. By a
general result
of Cheeger-Gromov-Taylor (see also Roe), the heat kernel $k^\pm(t,x,y)$
converges uniformly over compact subsets of $\widetilde M \times \widetilde M$
to $P^\pm(x,y)$, as $t \to \infty$.  Therefore one has
$$
\Ind_{1, \Gamma, \sigma}(\widetilde{\not\partial}^+_{\mathcal E}
\otimes \nabla)
= \tr (\theta_\sigma(R_t)) = \tr_s( e^{-t(
\widetilde{\not\partial}_{\mathcal E}
\otimes \nabla)^2})
$$
\begin{align*}
\lim_{t\to\infty} \tr(e^{-t(\widetilde{\not\partial}^\pm_{\mathcal E}
\otimes \nabla)^2})
& = \lim_{t\to\infty} \int_M \tr(k^\pm(t,x,x)) dx\\
& = \int_M \tr(P^\pm(x,x)) dx\\
& = \tr(P^\pm) \tag{1}
\end{align*}
Next observe that
\begin{align*}
\frac{\partial}{\partial t} \tr_s( e^{-t
(\widetilde{\not\partial}_{\mathcal E}
\otimes \nabla)^2})
& = - \tr_s( (\widetilde{\not\partial}_{\mathcal E}
\otimes \nabla)^2 e^{-t(\widetilde{\not\partial}_{\mathcal E}
\otimes \nabla)^2})\\
& = -  \tr_s( [\widetilde{\not\partial}_{\mathcal E}
\otimes \nabla, (\widetilde{\not\partial}_{\mathcal E}
\otimes \nabla) e^{-t (\widetilde{\not\partial}_{\mathcal E}
\otimes \nabla)^2}])\\
& = 0
\end{align*}
since $\widetilde{\not\partial}_{\mathcal E}\otimes \nabla$
is an odd operator.
Here $\tr_s$ denotes the graded trace, i.e. the composition of the trace
$\tr$ and the
grading operator. Therefore we deduce that
\begin{align*}
\tr_s( e^{-t (\widetilde{\not\partial}_{\mathcal E}
\otimes \nabla)^2}) & =
\lim_{t\to\infty} \tr_s( e^{-t (\widetilde{\not\partial}_{\mathcal E}
\otimes \nabla)^2})\\
& = \tr_s (P) \\
& = \Index_{L^2}(\widetilde{\not\partial}^+_{\mathcal E}
\otimes \nabla).
\tag{2}
\end{align*}
By the local index theorem of Atiyah-Bott-Patodi \cite{ABP},
Getzler \cite{Get}, one has
$$
\lim_{t\to 0}\left( \tr(k^+(t,x,x)) - \tr(k^-(t,x,x))   \right)
= [\hat{A}(\Omega)\,\tr(e^{R^{\E}}) e^{\omega}]_n\eqno{(3)}
$$
where $[\;\;]_n$ denotes the component of degree $n=$ dim $M$,
$\Omega$ is the curvature of the metric on $\widetilde M$,
$R^{\E}$ is the curvature of the connection on $\widetilde {\mathcal E}$.
Combining equations (1), (2) and (3), one has
\begin{align*}
\Ind_{1, \Gamma, \sigma}(\widetilde{\not\partial}^+_{\mathcal E}
\otimes \nabla)
& = \Index_{L^2}(\widetilde{\not\partial}^+_{\mathcal E}\otimes \nabla)\\
& = \int_M \hat{A}(\Omega)\,\tr(e^{R^{\E}}) e^{\omega}.
\end{align*}
We shall now generalize this argument to the case when
$c\in Z^{2q}(\Gamma)$, $q>0$. Here we adapt the strategy and
proof in \cite{CM} to our situation. By section 2.2, one has for $t>0$
\begin{align*}
\Ind_{c, \Gamma, \sigma}(\widetilde D\otimes \nabla)
& = \tr_c (\theta_\sigma(R_t),\theta_\sigma(R_t), \ldots \theta_\sigma(R_t))\\
& =\sum_{\gamma_0\gamma_1\ldots\gamma_{2q}=1}
\tr(\theta_\sigma(R_t)_{\gamma_0}\theta_\sigma(R_t)_{\gamma_1}
\ldots \theta_\sigma(R_t)_{\gamma_{2q}})
c(1, \gamma_1,\ldots,\gamma_1\gamma_2\ldots\gamma_{2q})
\tr(\delta_{\gamma_0}\ldots\delta_{\gamma_{2q}}).
\end{align*}
By changing variables, $\gamma_1 = \gamma_1, \gamma_2 = \gamma_1\gamma_2,
\ldots \gamma_{2q} = \gamma_1\ldots \gamma_{2q},$
one obtains
$$
\Ind_{c, \Gamma, \sigma}(\widetilde D\otimes \nabla)
$$
$$
=\sum_{\gamma_1,\gamma_{2}, \ldots \gamma_{2q}\in \Gamma}
\tr(\theta_\sigma(R_t)_{\gamma_1}\theta_\sigma(R_t)_{\gamma_1^{-1}\gamma_2}
\ldots\theta_\sigma(R_t)_{\gamma_{2q}^{-1}})
c(1, \gamma_1,\ldots,\gamma_{2q})
\tr(\delta_{\gamma_1}\delta{\gamma_1^{-1}\gamma_2}
\ldots\delta_{\gamma_{2q}^{-1}})
$$
$$
= \int_{M^{2q+1}} \sum_{\gamma_1,\ldots\gamma_{2q}\in \Gamma}
c(1, \gamma_1,\ldots,\gamma_{2q})
\tr(\delta_{\gamma_1}\delta{\gamma_1^{-1}\gamma_2}
\ldots\delta_{\gamma_{2q}^{-1}}) \times$$
$$
{\tr}(R_t(\beta(x_0), \gamma_1\beta(x_1)) \ldots
R_t(\gamma_{2q}\beta(x_{2q}), \beta(x_0)))dx_0\ldots dx_{2q}
$$
Observe that if $\beta_U : U \to \M$ is a smooth local
section, then there is a unique element $(\gamma_1, \ldots, \gamma_{2q})\in
\Gamma^{2q}$ such that
$(\beta(x_0), \gamma_1\beta(x_1),\ldots \gamma_{2q}\beta(x_{2q})) \in
\beta_U(U)^{2q+1}$, and in which case, one has the equality
$c(1, \gamma_1,\ldots,\gamma_{2q})  = \phi_c(x_0,x_1, \ldots x_{2q})$ \ \ (and
$\phi_c = 0$ otherwise),
where $\phi_c$ denotes the $\Gamma$-equivariant
(Alexander-Spanier) $2q$-cocycle on
$\widetilde M$ representing
the $f^*(c)$ which is the pullback of the group 2-cocycle
by the classifying map $f$,
where we have identified $M$ with a fundamental domain for the
$\Gamma$ action on $\widetilde M$, just as was done in the
case when the group cocycle $c = 1$. Since
$R_t$ is mainly supported near the diagonal as $t\to 0$,
and using the equivariance of $R_t$, one sees that
$$
\Ind_{c, \Gamma, \sigma}(\widetilde D\otimes \nabla)
= \int_{M^{2q+1}} \phi_c(x_0,x_1, \ldots x_{2q})
{\tr}(R_t(x_0, x_1)\ldots R_t(x_{2q}, x_0))dx_0dx_1\ldots dx_{2q}
$$
The proof is completed by taking the limit as
$t\to 0$ and by applying the {\em local} higher index
Theorems 3.7 and 3.9 in \cite{CM}.
\end{proof}

\begin{rems*}
A particular case of Theorem 2.6 highlights a key new phenomenon in
the case of orbifolds, viz. in the special case when the group cocycle
$c =1 \in Z^0(\Gamma)$ is trivial, and when the multiplier $\sigma=1$
is trivial, then $\Ind_{1, \Gamma, 1}(\widetilde D)$
is the $L^2$ index of $\widetilde D$ as defined Atiyah.
By comparing with the cohomological formula due to Kawasaki \cite{Kaw}
for the Fredholm index of the operator $ D$ on the orbifold $M$, we see
that in general these are {\em not} equal, and the error term is a
rational number which can be
expressed explicitly as a cohomological formula on the lower dimensional
strata
of the orbifold $M$. Since Atiyah's $L^2$ index theorem can be viewed
as an integrality statement for the $L^2$ index of $\widetilde D$,
$\Ind_{1, \Gamma, 1}(\widetilde D)$ in the smooth
case, we therefore see that for general orbifolds the $L^2$ index of
$\widetilde D$,
$\Ind_{1, \Gamma, 1}(\widetilde D)$ is only a {\em rational} number.
This was also observed by \cite{Far}.
Whereas for smooth two dimensional manifolds
the {\em higher} index associated to the area
cocycle is an {\em integer} \cite{CHMM},
in section 5, we show that it is however only a {\em
rational} number for general two dimensional orbifolds.
In some other work in progress, we will
give an alternate heat kernel proof
of a generalization of this theorem, using superconnections.
\end{rems*}

\section{Range of the trace and the Kadison constant}

In this section, we will first calculate the range of the
canonical trace map on $K_0$ of the twisted group $C^*$-algebras
for Fuchsian groups $\Gamma$ of signature $(g,\nu_1,\ldots ,\nu_n)$.
We use in an essential way some of the results of the previous section
such as the twisted version of the $L^2$-index theorem of Atiyah
\cite{At}, which is due to Gromov \cite{Gr2}, and which is proved
in Theorem 2.6.
This enables us to deduce information about projections in the
twisted group $C^*$-algebras. In the case of no twisting, this follows
because the Baum-Connes conjecture is known to be true while
these results are also well
known for the case of the irrational rotation algebras, and
for the twisted groups $C^*$ algebras of the
fundamental groups of closed Riemann surfaces of positive genus
\cite{CHMM}. Our theorem generalises most of these results.
We will apply the results of this section in the next section to study
some quantitative aspects of the
spectrum of projectively periodic elliptic operators,
mainly on orbifold covering spaces
of hyperbolic orbifolds.

\subsection{The isomorphism classes of algebras ${C}^*(\Gamma, \sigma)$}
Let $\sigma \in Z^2(\Gamma, U(1))$ be a multiplier on $\Gamma$,
where $\Gamma$ is a Fuchsian group of signature $(g, \nu_1, \ldots, \nu_n)$. If
$\sigma' \in Z^2(\Gamma, U(1))$ is another multiplier on $\Gamma$ such
that $[\sigma] = [\sigma'] \in H^2(\Gamma, U(1))$, then it can be easily
shown that ${C}^*(\Gamma, \sigma) \cong {C}^*(\Gamma, \sigma')$. That is,
the isomorphism classes of the $C^*$-algebras ${C}^*(\Gamma_{g'}, \sigma)$ are
naturally parametrized by $H^2(\Gamma, U(1))$. In particular, if we consider
only multipliers $\sigma$ such that $\delta(\sigma) = 0$, we see that
these are parametrised by $\text{ker} (\delta) \subset H^2(\Gamma, U(1))$.
It follows from the discussion at the beginning of the next subsection
that $\text{ker} (\delta) \cong U(1)$.
We summarize this below.

\begin{lemma} Let $\Gamma$ be a Fuchsian group of signature
$(g, \nu_1, \ldots, \nu_n)$.
Then the isomorphism classes of twisted group $C^*$-algebras
${C}^*(\Gamma, \sigma)$ such that $\delta(\sigma) = 0$ are
naturally parametrized by $U(1)$.
\end{lemma}

\subsection{K-theory of twisted group $C^*$ algebras}
We begin by computing the $K$-theory of twisted group $C^*$-algebras
for Fuchsian groups $\Gamma$ of signature $(g,\nu_1,\ldots ,\nu_n)$.
Let $\sigma$ be a multiplier on $\Gamma$. It defines a cohomology class
$[\sigma] \in H^2(\Gamma, U(1))$.
Consider now the short
exact sequence of coefficient groups
\[
   1\to\mathbb{Z} \overset{i}{\to} \mathbb{R}
      \overset{e^{2\pi\sqrt{-1}}}{\longrightarrow} U(1) \to 1,
\]
which gives rise to a long exact sequence of cohomology groups (the change of
coefficient groups sequence)
$$
   \cdots \to H^2(\Gamma,\mathbb{Z}) \overset{i_*}{\to}
H^2(\Gamma,\mathbb{R})
      \overset{{e^{2\pi\sqrt{-1}}}_*}{\longrightarrow} H^2(\Gamma, U(1))
      \overset{\delta}{\to} H^3(\Gamma,\mathbb{Z})\overset{i_*}{\to}
H^3(\Gamma,\mathbb{R}).\eqno{(3.1)}
$$

We first show that the map
$$ H^2(\Gamma, U(1))\overset{\delta}{\to} H^3(\Gamma,\mathbb{Z})$$
is a a surjection.

In fact, it is enought to show that $H^3(\Gamma,\mathbb{R})=\{0\}$. In
order to see this it is enough to notice that we have
a $G$ action on $B\Gamma_{g'}$ with quotient $B\Gamma$,
$$ G\to B\Gamma_{g'} \stackrel{\lambda}{\to} B\Gamma \eqno{(3.2)}$$
and therefore the Leray-Serre spectral sequence, we have
$$ E^2=Tor^{H_*(G, \R)}(\R,H_*(B\Gamma_{g'},\R)) $$
that converges to $H_*(B\Gamma,\R)$. Moreover, we have
$$ E^2=Tor^{H_*(G, \R)}(\R,\R) $$
converging to $H_*(BG,\R)$, \cite{McCl} 7.16.
Notice also that, with $\R$ coefficients, we have $
H_q(BG,\R)=\{0\}$ for
$q>0$. Thus we obtain that, with $\R$ coefficients,
$H_q(B\Gamma,\R)\cong H^q(B\Gamma, \R)$
is $\R$ in degrees $q=0$ and $q=2$, $\R^{2g}$ in degree $q=1$, and
trivial in degrees $q>2$. I particular, $(3.1)$ now becomes
$$
   \cdots \to H^2(\Gamma,\mathbb{Z}) \overset{i_*}{\to}
H^2(\Gamma,\mathbb{R})
      \overset{{e^{2\pi\sqrt{-1}}}_*}{\longrightarrow} H^2(\Gamma, U(1))
      \overset{\delta}{\to} H^3(\Gamma,\mathbb{Z})\overset{i_*}{\to}
0.\eqno{(3.3)}
$$

In the following $[\Gamma]$ will denote a choice of a generator in
$H_2(B\Gamma,\R)\cong\R\cong H^2( B\Gamma, \R)$. Using equation $(3.2)$
and the previous argument, we see that $\lambda_*[\Sigma_{g'}]
= \#(G) [\Gamma]$, since $B\Gamma_{g'}$ and $\Sigma_{g'}$ are homotopy
equivalent,
and where $ \#(G)$ denotes the order of the finite group $G$.

In particular, for any multiplier $\sigma$ of $\Gamma$ with $ [\sigma]
\in H^2(\Gamma_{g'}, U(1))$ and with $\delta(\sigma) =0$,
there is a $\R$-valued 2-cocycle $\zeta$ on $\Gamma$ with
$[\zeta] \in H^2(\Gamma, \R)$ such that
$[e^{2\pi\sqrt{-1}\zeta}] = [\sigma]$. Define a homotopy
$[\sigma_t] = [e^{2\pi\sqrt{-1}t\zeta}] \quad \forall t \in [0,1]$ which
is a homotopy of multipliers $\sigma_t$ that connects the multiplier
$\sigma$ and the trivial multiplier. Note also that this homotopy is {\em
canonical}
and not dependent on the particular choice of $\zeta$. Therefore one obtains a
homotopy of ismorphism classes of twisted group $C^*$-algebras
$C^*(\Gamma, \sigma_t)$ connecting $C^*(\Gamma, \sigma)$ and
$C^*(\Gamma)$. It is this homotopy which will essentially be
used to show that $C^*(\Gamma, \sigma)$ and $C^*(\Gamma)$ have the
same $K$-theory.

Let $\Gamma\subset G$ be a discrete
cocompact subgroup of $G$ and $A$ be an algebra
admitting an action of $\Gamma$ by automorphisms.
Then the cross product algebra $[A\otimes C_0(G)]\rtimes \Gamma$,
is Morita equivalent to the algebra of continuous sections
vanishing at infinity
$C_0(\Gamma\backslash G, \mathcal{E})$, where $\mathcal{E}\to
\Gamma\backslash G$ is the flat $A$-bundle defined as the quotient
\begin{equation*}
   \mathcal{E} = (A\times G)/\Gamma \to \Gamma\backslash G.
\end{equation*}
Here we consider the diagonal action of $\Gamma$ on $A\times G$.
We refer the reader to \cite{Kas1} for the technical defintion of
a $K$-amenable group. However we mention
that any solvable Lie group, and in fact any amenable Lie group
is $K$-amenable, and in fact it is shown in
\cite{Kas1}, \cite{JuKas} that the
non-amenable groups ${\mathbf{SO}}_0(n,1)$, ${\mathbf{SU}}(n,1)$
are $K$-amenable Lie groups. Also,
Cuntz \cite{Cu} has shown that the class of $K$-amenable groups
is closed under the operations of taking subgroups, under free products
and under direct products.

\begin{thm}[\cite{Kas1},\cite{Kas2}]
If $G$ is $K$-amenable, then
$(A\rtimes\Gamma)\otimes C_0(G)$ and $[A\otimes C_0(G)]\rtimes
\Gamma$
have the same $K$-equivariant $K$-theory, where $K$ acts in the standard
way on $G$ and trivially on the other factors.
\end{thm}

Combining Theorem 3.2 with the remarks above, one gets the following
important corollary.

\begin{cor}
If $G$ is $K$-amenable, then $(A\rtimes\Gamma)\otimes C_0(G)$ and
$C_0(\Gamma\backslash G, \mathcal{E})$ have the same
$K$-equivariant $K$-theory.
Equivalently, one has for $j=0,1$,
$$
{K_K}_j(C_0(\Gamma\backslash G, \mathcal{E})) \cong {K_K}_{j+ \dim(G/K)}
(A\rtimes\Gamma).
$$
\end{cor}

We now come to the main theorem of this section,
which generalizes theorems of \cite{CHMM}, \cite{PR}, \cite{PR2}.

\begin{thm}
Suppose that $\Gamma$ is a discrete cocompact subgroup
in a $K$-amenable Lie group $G$  and that
$K$ is a maximal compact subgroup of $G$. Then
$$
{K}_\bullet(C^*(\Gamma,\sigma)) \cong
{K_K}^{\bullet + \dim(G/K)}(\Gamma\backslash G, \delta(B_\sigma)),
$$
where $\sigma \in H^2(\Gamma, U(1))$ is any multiplier on $\Gamma$,
${K_K}^{\bullet}(\Gamma\backslash G,
\delta(B_\sigma))$ is the twisted $K$-equivariant
$K$-theory of a continuous trace
 $C^*$-algebra
$B_\sigma$ with spectrum $\Gamma\backslash G$, and $\delta(B_\sigma)$
denotes the
Dixmier-Douady invariant of $B_\sigma$.\end{thm}

\begin{proof}
Let $\sigma\in H^2(\Gamma, U(1))$, then the {\em twisted} cross product
algebra $A\rtimes_\sigma \Gamma$ is stably equivalent to the
cross product $(A\otimes
\mathcal{K})\rtimes\Gamma$ where $\mathcal{K}$ denotes compact operators.
This is the Packer-Raeburn stabilization trick \cite{PR}, which we now
describe in more detail.
Let $V:
\Gamma \to U(\ell^2(\Gamma))$ denote the left regular $(\Gamma, \bar\sigma)$
representation on $\ell^2(\Gamma)$, i.e. for $\gamma, \gamma_1 \in \Gamma$
and $f\in \ell^2(\Gamma)$
$$
(V(\gamma_1)f)(\gamma) = \bar\sigma(\gamma_1, \gamma_1^{-1}\gamma)
f(\gamma_1^{-1}\gamma).
$$
Then for $\gamma_1, \gamma_2 \in \ell^2(\Gamma)$,  $V$
satisfies $V(\gamma_1)V(\gamma_2) = \bar\sigma(\gamma_1,
\gamma_2)V(\gamma_1\gamma_2)$. That is, $V$ is a projective representation
of $\Gamma$. Since $Ad$ is trivial on $U(1)$, it follows that
$\alpha(\gamma) = Ad(V(\gamma))$ is a representation of $\Gamma$
on $\mathcal K$. This is easily generalised to the case when
$\mathbb C$ is replaced by the $*$ algebra $A$.

Using Corollary 3.3 again,
one sees
that $A\rtimes_\sigma\Gamma \otimes C_0(G)$ and $C_0(\Gamma\backslash
G, \mathcal{E}_\sigma)$
have the same $K$-equivariant $K$-theory, whenever $G$ is $K$-amenable, where
$$
   \mathcal{E}_\sigma = (A\otimes\mathcal{K} \times G)/\Gamma \to
      \Gamma\backslash G
$$
is a flat $A\otimes\mathcal{K}$-bundle over $\Gamma\backslash G$ and
$K$ is a maximal compact subgroup of $G$.
In the particular case when $A= \mathbb{C}$, one sees that $C^*_r(\Gamma,
\sigma)\otimes C_0(G)$ and $C_0(\Gamma\backslash G,
\mathcal{E}_\sigma)$ have the same
$K$-equivariant $K$-theory whenever $G$ is $K$-amenable, where
\[
   \mathcal{E}_\sigma = (\mathcal{K} \times G)/\Gamma\to\Gamma\backslash G.
\]
But the twisted $K$-equivariant
$K$-theory ${K_K}^{*}(\Gamma\backslash G, \delta(B_\sigma))$ is
by definition the  $K$-equivariant $K$-theory of the
continuous trace $C^*$-algebra $B_\sigma
= C_0(\Gamma\backslash G, \mathcal{E}_\sigma)$ with spectrum
$\Gamma\backslash G$.
Therefore
$$
K_\bullet(C^*(\Gamma,\sigma)) \cong {K_K}^{\bullet + \dim(G/K)}
(\Gamma\backslash G, \delta(B_\sigma)).
$$

\end{proof}

Our next main result says that for discrete cocompact
subgroups in $K$-amenable Lie groups,
the reduced and unreduced twisted group $C^*$-algebras have canonically
isomorphic $K$-theories. Therefore all the results that we prove regarding
the $K$-theory of these reduced twisted group $C^*$-algebras are also valid
for the unreduced twisted group $C^*$-algebras.

\begin{thm}
Let $\sigma\in H^2(\Gamma, U(1))$ be a multiplier on $\Gamma$ and
$\Gamma$ be a discrete cocompact
subgroup in a $K$-amenable Lie group. Then the canonical
morphism $C^*(\Gamma, \sigma) \rightarrow C^*_r(\Gamma, \sigma)$ induces
an isomorphism
$$
K_*(C^*(\Gamma, \sigma)) \cong K_*(C^*_r(\Gamma, \sigma)).
$$
\end{thm}

\begin{proof}
We note that by the Packer-Raeburn trick, one has
$$
C^*(\Gamma, \sigma) \otimes {\mathcal K} \cong  {\mathcal K}  \rtimes  \Gamma
$$
and
$$
C^*_r(\Gamma, \sigma) \otimes {\mathcal K} \cong  {\mathcal K}  \rtimes_r
\Gamma,
$$
where $\rtimes_r$ denotes the reduced crossed product.
Since $\Gamma$ is a lattice in a $K$-amenable Lie group,
the canonical
morphism ${\mathcal K}  \rtimes  \Gamma \rightarrow
{\mathcal K}  \rtimes_r  \Gamma$ induces
an isomorphism (cf. \cite{Cu})
$$
K_*({\mathcal K}  \rtimes  \Gamma) \cong K_*({\mathcal K}  \rtimes_r  \Gamma),
$$
which proves the result.
\end{proof}

We now specialize to the case when $G= \so_0(2,1)$, $K=\so(2)$ and
$\Gamma = \Gamma(g, \nu_1, \ldots, \nu_n)$ is a Fuchsian
group, i.e. the orbifold
fundamental group of a hyperbolic orbifold of signature
$(g, \nu_1, \ldots, \nu_n)$, \ $\Sigma(g, \nu_1, \ldots, \nu_n)$,
where  $\Gamma \subset G$ (note that $G$ is $K$-amenable), or
when $G= \mathbb R^2$, $K= \{e\}$ and $g=1$, with $\Gamma$ being
a cocompact crytallographic group.

\begin{prop}
Let $\sigma$ be a multiplier on the Fuchsian group $\Gamma$ of signature
$(g,\nu_1,\ldots ,\nu_n)$ such that $\delta(\sigma) = 0$.
Then one has
\begin{enumerate}
\item $K_0(C^*(\Gamma, \sigma)) \cong K_0(C^*(\Gamma)) \cong
K^0_{orb}(\Sigma(g, \nu_1,\ldots,\nu_n)) \cong \Z^{2-n +\sum_{j=1}^n \nu_j}$
\item $K_1(C^*(\Gamma, \sigma)) \cong K_1 (C^*(\Gamma)) \cong
K^1_{orb}(\Sigma(g, \nu_1,\ldots,\nu_n)) \cong \mathbb{Z}^{2g}$.
\end{enumerate}

\end{prop}

\begin{proof}
Now by a result due to Kasparov \cite{Kas1}, which he proves by connecting the
regular representation to the trivial one via the complementary series,
one has
$$
K_\bullet (C^*(\Gamma)) \cong K^\bullet_{SO(2)}(P(g, \nu_1,\ldots,\nu_n)) =
K^\bullet_{orb}(\Sigma(g,\nu_1,\ldots ,\nu_n)).
$$
We recall next the calculation of Farsi \cite{Far}
for the orbifold
$K$-theory of the hyperbolic orbifold $\Sigma(g,\nu_1,\ldots ,\nu_n)$
$$
K^0_{orb}(\Sigma(g, \nu_1,\ldots,\nu_n))\equiv K_0(C^*(\Sigma(g,
\nu_1,\ldots,\nu_n)) ) =
K^0_{SO(2)}(P(g, \nu_1,\ldots,\nu_n))
\cong \Z^{2-n +\sum_{j=1}^n \nu_j}
$$
and
$$
K^1_{orb}(\Sigma(g, \nu_1,\ldots,\nu_n))\equiv K_1(C^*(\Sigma(g,
\nu_1,\ldots,\nu_n)) ) =
K^1_{SO(2)}(P(g, \nu_1,\ldots,\nu_n))
\cong \Z^{2g}
$$

By Theorem 3.4 we have
\[
   K_j(C^*(\Gamma)) \cong K^j_{SO(2)}(P(g, \nu_1,\ldots,\nu_n))
\quad\text{for }j=0,1,
\]
and more generally
\[
   K_j(C^*(\Gamma,\sigma)) \cong K^j_{SO(2)}(P(g, \nu_1,\ldots,\nu_n),
\delta(B_\sigma)),\quad
      j = 0,1,
\]
where $B_\sigma = C(P(g, \nu_1,\ldots,\nu_n), \mathcal{E}_\sigma)$.
Finally, because $\mathcal{E}_\sigma$ is a locally trivial bundle of
 $C^*$-algebras over $P(g, \nu_1,\ldots,\nu_n)$,
with fibre $\mathcal{K}$ ($=$ compact operators),
it has a Dixmier-Douady invariant $\delta(B_\sigma)$ which can be
viewed as the obstruction to $B_\sigma$ being Morita
equivalent to $C(\Sigma_g)$.  But by assumption
$\delta(B_\sigma) = \delta(\sigma) = 0$.
Therefore $B_\sigma$ is Morita equivalent to $C(P(g, \nu_1,\ldots,\nu_n))$
and we conclude that
\[
   K_j(C^*(\Gamma, \sigma)) \cong  K^j_{SO(2)}(P(g, \nu_1,\ldots,\nu_n))
\cong K^j_{orb}(\Sigma(g, \nu_1,\ldots,\nu_n))\quad j = 0,1.
\]
\end{proof}

\subsection{Twisted Kasparov map}
Let $\Gamma$ be as before,
that is, $\Gamma$ is the orbifold fundamental group of
the hyperbolic orbifold $\Sigma(g, \nu_1,\ldots,\nu_n)$.
Then for any multiplier $\sigma$ on $\Gamma$, the
\emph{twisted Kasparov isomorphism},
$$
 \mu_\sigma : K^\bullet_{orb} (\Sigma(g, \nu_1,\ldots,\nu_n)) \to
K_\bullet (C^*_r(\Gamma,\sigma))\eqno{(3.4)}
$$
is defined as follows.  Let
$\mathcal{E}\to\Sigma(g, \nu_1,\ldots,\nu_n)$ be an orbifold
vector bundle over $\Sigma(g, \nu_1,\ldots,\nu_n)$
defining  an element $[\mathcal{E}]$
in $K^0(\Sigma(g, \nu_1,\ldots,\nu_n))$. As in \cite{Kaw}, one can form the
twisted Dirac
operator $\npartial^+_{\mathcal{E}} : L^2(\Sigma(g, \nu_1,\ldots,\nu_n),
\mathcal{S}^+\otimes\mathcal{E})
\to L^2(\Sigma(g, \nu_1,\ldots,\nu_n), \mathcal{S}^-\otimes\mathcal{E})$
where $\mathcal{S}^{\pm}$ denote the $\frac12$
spinor bundles over $\Sigma(g, \nu_1,\ldots,\nu_n)$.
By Proposition 3.2 of the previous subsection,
there is a canonical isomorphism
$$
K_\bullet (C^*_r(\Gamma, \sigma)) \cong K^\bullet_{orb} (\Sigma(g,
\nu_1,\ldots,\nu_n)).
$$
Both of these maps are assembled to yield the twisted Kasparov map as in
$(3.4)$.
Observe that $\Sigma(g, \nu_1,\ldots,\nu_n) = \underline B\Gamma$, and that
the twisted Kasparov map has a natural generalisation, which will
be studied elsewhere.

We next describe this map more explicitly.
One can lift the twisted Dirac operator $\npartial^+_{\mathcal{E}}$ as above,
to a $\Gamma$-invariant operator $\widetilde{\npartial_{\mathcal{E}}^+}$
on  $\mathbb{H} = \widetilde\Sigma(g, \nu_1,\ldots,\nu_n)$, which is
the universal orbifold cover of $\Sigma(g, \nu_1,\ldots,\nu_n)$,
\[
   \widetilde{\npartial_{\mathcal{E}}^+} : L^2(\mathbb{H},
      \widetilde{\mathcal{S}^+\otimes\mathcal{E}}) \to L^2(\mathbb{H},
      \widetilde{\mathcal{S}^-\otimes\mathcal{E}})
\]
Therefore as before in $(3.3)$,
for any multiplier $\sigma$ of $\Gamma$ with $\delta([\sigma]) =1$,
there is a $\R$-valued 2-cocycle $\zeta$ on $\Gamma$ with
$[\zeta] \in H^2(\Gamma, \R)$ such that
$[e^{2\pi\sqrt{-1}\zeta}] = [\sigma]$. By the earlier argument using spectral
sequences and the fibration as in equation $(3.2)$, we see that the map $
\lambda$ induces an isomorphism
$H^2(\Gamma, \R) \cong H^2(\Gamma_{g'}, \R)$, and therefore
there is a 2-form $\omega$ on $\Sigma_{g'}$ such that
$[e^{2\pi\sqrt{-1}\omega}] = [\sigma]$.
Of course, the choice of $\omega$ is not unique, but this will not affect the
results that we are concerned with. Let $\widetilde \omega$ denote the lift
of $\omega$ to the universal cover $\mathbb H$. Since the hyperbolic plane
$\mathbb H$ is contractible, it follows that $\widetilde \omega = d\eta$
where $\eta$ is a 1-form on $\mathbb H$ which is not in general $\Gamma$
invariant. Now $\nabla = d - i\eta$ is a Hermitian connection on the trivial
complex line bundle on $\mathbb H$. Note that the curvature of $\nabla$ is
$\nabla^2 = i\tilde\omega$. Consider now the twisted Dirac operator
$\widetilde \npartial^+_{\mathcal{E}}$ which is
twisted again by the connection $\nabla$,
\[
   \widetilde{\npartial_{\mathcal{E}}^+}\otimes\nabla : L^2(\mathbb{H},
      \widetilde{\mathcal{S}^+\otimes\mathcal{E}}) \to L^2(\mathbb{H},
      \widetilde{\mathcal{S}^-\otimes\mathcal{E}}).
\]
It does not commute with the $\Gamma$ action, but it does commute
with the projective $(\Gamma, \sigma)$-action  which is defined by the
multiplier
$\sigma$, and by the twisted $L^2$-index theorem of the previous section,
it has a
$\Gamma$-$L^2$-index
\[
   \Ind_{(\Gamma, \sigma)} (\widetilde{\npartial^+_{\mathcal{E}}}\otimes
\nabla) \in
      K_0(\mathcal{R}(\Gamma, \sigma)),
\]
where  as before, $\mathcal{R}$ denotes the algebra of rapidly decreasing
sequences
on $\Z^2$.
Then observe that the {\em twisted Kasparov map} is
$$
\mu_\sigma([\mathcal{E}]) = j_*(\Ind_{(\Gamma, \sigma)}
(\widetilde{\npartial^+_{\mathcal{E}}}\otimes \nabla) ) \in
K_0({C}^*(\Gamma, \sigma)),
$$
where $j: {\mathcal R}(\Gamma, \sigma) = {\C}(\Gamma, \sigma) \otimes
\mathcal R  \to
{C}^*_r(\Gamma, \sigma)\otimes \mathcal K$
is the natural inclusion map, and
$$j_* : K_0 ( {\mathcal R}(\Gamma, \sigma)) \to K_0 ({C}^*_r(\Gamma,
\sigma))$$ is the
induced map on $K_0$.

The canonical trace on ${C}^*_r(\Gamma, \sigma))$ induces a linear map
$$
[\tr] : K_0 ({C}^*_r(\Gamma, \sigma)) \to \mathbb R
$$
which is called the {\em trace map} in $K$-theory.
Explicitly, first $\tr$ extends to matrices with entries in
${C}^*(\Gamma, \sigma)$ as (with Trace denoting matrix trace):
\[
   \tr(f\otimes r) = {\mbox{Trace}}(r) \tr(f).
\]

Then the extension of $\tr$ to $K_0$ is given by
$[\tr]([e]-[f]) = \tr(e) - \tr(f)$, where $e,f$ are idempotent matrices with
entries in ${C}^*(\Gamma, \sigma))$.

\subsection{Range of the trace map on $K_0$}
We can now state the major theorem of this section.

\begin{thm} Let $\Gamma$ be a Fuchsian group of signature $(g, \nu_1, \ldots
, \nu_n)$, and $\sigma$ be a multiplier of $\Gamma$ such that
$\delta(\sigma) = 0$.
Then the range of the trace map is
$$
[\tr] (K_0 ({C}^*_r(\Gamma, \sigma)) ) = \Z\theta + \Z + \sum_{i=1}^n \Z
(1/\nu_i) ,
$$
where $2\pi\theta = \langle[\sigma], [\Gamma]\rangle\ \in (0,1]$ is the
result of
pairing the
multiplier $\sigma$ with the fundamental class of
$\Gamma$ (cf. subsection 3.1).
\end{thm}

\begin{proof}

We first observe that by the results of the previous subsection the
twisted Kasparov map is an isomorphism. Therefore to compute the range of the
trace map on $K_0$, it suffices to compute the range of the trace map
on elements of the form
$$\mu_\sigma([\mathcal{E}^0 ] -
[\mathcal{E}^1])$$
 for any element
 $$[\mathcal{E}^0 ] -
[\mathcal{E}^1] \in K^0_{orb}(\Sigma(g, \nu_1,\ldots,\nu_n)).$$
where $\mathcal{E}^0, \mathcal{E}^1$ are orbifold vector bundles over
the orbifold $\Sigma(g, \nu_1,\ldots,\nu_n)$, which as in section 1, can be
viewed
as $G$-equivariant vector bundles over the Riemann surface $\Sigma_{g'}$
which is
an orbifold $G$ covering of the orbifold $\Sigma(g, \nu_1,\ldots,\nu_n)$.

By the twisted  $L^2$ index theorem for orbifolds from the previous section,
one has
$$
   [\tr](\Ind_{(\Gamma, \sigma)} (\widetilde{\npartial_{\mathcal{E}}^+}\otimes
\nabla)) =
       \frac{1}{2\pi}\int_{\Sigma(g, \nu_1,\ldots,\nu_n)} \hat{A}(\Omega)
\tr(e^{R^{\mathcal{E}}}) e^{\omega}.
      \eqno{(3.5)}
$$
We next simplify the right hand side of equation $(3.5)$ using
\begin{align*}
   \hat{A}(\Omega) &= 1  \\
   \tr(e^{R^{\mathcal{E}}}) &= \rank \mathcal{E} + \tr(R^{\mathcal{E}}) \\
   e^{\omega} &= 1 + {\omega} .
\end{align*}
Therefore one has
\[
   [\tr](\Ind_{(\Gamma, \sigma)} (\widetilde{\npartial_{\mathcal{E}}^+}\otimes
\nabla)) =
      {\rank \mathcal E}\frac{\int_{\Sigma(g, \nu_1,\ldots,\nu_n)}
\omega}{2\pi}
      + \frac{ \int_{\Sigma(g, \nu_1,\ldots,\nu_n)}\tr(R^{\mathcal{E}})}{2\pi},
\]
Now by the index theorem for orbifolds, due to Kawasaki \cite{Kaw}, we see
that
$$
 \frac{ \int_{\Sigma(g, \nu_1,\ldots,\nu_n)}\tr(R^{\mathcal{E}})}{2\pi}
+ \frac{\sum_{i=1}^n \beta_i/\nu_i }{2\pi}  = \text{index}
(\npartial_{\mathcal{E}}^+)
\in \Z,
$$
Therefore we see that
$$
 \frac{ \int_{\Sigma(g, \nu_1,\ldots,\nu_n)}\tr(R^{\mathcal{E}})}{2\pi}
\in \Z + \sum_{i=1}^n \Z (1/\nu_i)
$$
Observe that
$$
\int_{\Sigma(g, \nu_1,\ldots,\nu_n)} \omega
= \frac{1}{\#(G) }\int_{\Sigma_{g'}} \omega = \langle [\omega], [\Gamma]
\rangle
$$ since $\Sigma_{g'}$ is
an orbifold $G$ covering of the orbifold $\Sigma(g, \nu_1,\ldots,\nu_n)$
and $[\Gamma]$ is equal to $\frac{[\Sigma_{g'}]}{\#(G)}$, cf. section 3.1 and
that by assumption,
$$\frac{\langle [\omega], [\Gamma] \rangle}{2\pi} -
\theta\in \mathbb Z.$$
It follows that the range of the trace map on $K_0$ is
$\mathbb Z \frac{\langle [\omega],
[\Gamma] \rangle}{2\pi} + \Z + \sum_{i=1}^n \Z (1/\nu_i)=
\mathbb Z \theta + \Z + \sum_{i=1}^n \Z (1/\nu_i)$.
\end{proof}

We will now discuss one application of this result,
leaving further applications to the
next section.The application studies the number of
projections in the twisted group $C^*$-algebra,
which is a problem of independent interest.

\begin{prop}
Let $\sigma$ be a multiplier on
$\Gamma$ such that $\delta(\sigma) = 0$,  and
$2\pi\theta= \langle[\sigma], [\Gamma]\rangle \in (0,1]$
be the result of pairing the cohomology class of
$\sigma$ with the fundamental class of $\Gamma$. If $\theta$ is
rational, then
there are at most a finite number of unitary equivalence classes of
projections, other than $0$ and $1$,
in the reduced twisted group $C^*$-algebra ${C}^*_r(\Gamma, \sigma)$.
\end{prop}

\begin{proof} By assumption, $\theta = p/q$.
Let $P$ be a projection in ${C}^*_r(\Gamma, \sigma)$. Then
$1-P$ is also a projection in ${C}^*_r(\Gamma, \sigma)$ and one has
$$
1= \tr(1) = \tr(P) + \tr(1-P).
$$
Each term in the above equation is non-negative. Since $\sigma$ is rational
and by
Theorem 3.3, it follows that the Kadison constant $C_\sigma(\Gamma)>0$
(see section 4 for the definition) and
$\tr(P) \in \{0, C_\sigma(\Gamma), 2C_\sigma(\Gamma), \ldots 1\}$.
By faithfulness and normality  of the
trace $\tr$, it follows that there are at most a finite number of
unitary equivalence classes of projections, other than those of
$0$ and $1$
in ${C}^*_r(\Gamma, \sigma)$.
\end{proof}

\section{Applications to the spectral theory of projectively periodic
elliptic operators and the classification of twisted group $C^*$ algebras}

In this section, we apply the range of the trace theorem of section 3, to
prove some qualitative results on the spectrum of projectively periodic
self adjoint elliptic operators on the universal covering
of a good orbifold. In particular,
we study generalizations of the hyperbolic
analogues of the Ten Martini Problem in \cite{CHMM} and
the Bethe-Sommerfeld conjecture.
We also classify upto isomorphism,
the twisted group $C^*$ algebras for a cocompact
Fuchsian group.

Let $M$ be a compact, good orbifold, that is, the universal cover
$\Gamma \to \M \to M$ is a smooth manifold
and we will assume as before that there is a $(\Gamma, \bar\sigma)$-action
on $L^2(\M)$ given by $T_\gamma = U_\gamma \circ S_\gamma \, \forall \gamma
\in \Gamma$. Let $\E, \ \F$ be Hermitian vector bundles on $M$ and let
$\E, \ \F$ be the corresponding lifts to $\Gamma$-invariants Hermitian vector
bundles on $\M$. Then there are $(\Gamma, \sigma)$-actions
on $L^2(\M, \E)$ and $L^2(\M, \F)$ which are also given by
$T_\gamma = U_\gamma \circ S_\gamma \, \forall \gamma \in \Gamma$.

Now let $D : L^2(\M, \E) \to L^2(\M, \F)$ be a self adjoint elliptic
differential
operator that commutes with the $(\Gamma, \bar\sigma)$-action
that was defined earlier.
We begin with some basic facts about the spectrum of
such an operator. Recall that the {\em discrete spectrum}
of $D$, $spec_{disc}(D)$ consists of all the eigenvalues of $D$ that have
finite multiplicity, and the {\em essential spectrum} of
$D$, $spec_{ess}(D)$ consists of the complement $spec(D) \setminus
spec_{disc}(D)$. That is, $spec_{ess}(D)$ consists of the
set of accumulation points of the spectrum of $D$, $spec(D)$.
Our first goal is to prove that the essential spectrum is
unbounded. Our proof will be a modification of an argument in \cite{Sar}.

\begin{lemma}
Let $D: L^2(\M, \E) \to L^2(\M, \E)$ be a self adjoint elliptic differential
operator that commutes with the $(\Gamma, \bar\sigma)$-action.
Then the {\em discrete spectrum}
of $D$ is empty.
\end{lemma}

\begin{proof}
Let $\lambda$ be an eigenvalue of $D$ nd $V$ denote the
corresponding eigenspace. Then $V$ is a $(\Gamma, \sigma)$-
invariant subspace of $L^2(\M, \E)$.
If $\mathcal F$ is a fundamental domain for the action of $\Gamma$ on $\M$,
one sees that
$$
L^2(\M, \E) \cong L^2(\Gamma)\otimes L^2(\mathcal F, \E|_{\mathcal F})
$$
which can be proved by choosing a bounded measurable almost everywhere smooth
section of the orbifold covering
$\M \to M$. Here
$(\Gamma, \bar\sigma)$-action on  $L^2(\mathcal F, \E|_{\mathcal F})$
is trivial, and is the left regular
$(\Gamma, \bar\sigma)$-action on $L^2(\Gamma)$.
Therefore it suffices to show that
the dimension of {\em any}  $(\Gamma, \bar\sigma)$-
invariant subspace $V$ of $L^2(\Gamma)$ is infinite dimensional.
Let $\{v_1, \ldots, v_N\}$ be
an orthonormal basis for $V$. Then one has
$$
T_\gamma v_i(\gamma')
= \sum_{j=1}^N U_{ij}(\gamma) v_j(\gamma') \qquad \forall
\gamma, \gamma' \in \Gamma
$$
where $U = (U_{ij}(\gamma) )$ is some $N\times N$ unitary
matrix. Therefore
\begin{align*}
N= \sum_{j=1}^N ||v_i||^2 & = \sum_{j=1}^N
\sum_{\gamma\in \Gamma} |v_i(\gamma\gamma')|^2\\
&=  \sum_{\gamma\in \Gamma} \sum_{i=1}^N
\sum_{j=1}^N \sum_{k=1}^N
U_{ij}(\gamma) \overline{U_{ik}(\gamma)}
v_j(\gamma')\overline{v_k(\gamma')}\\
& = \sum_{\gamma\in \Gamma} \sum_{j=1}^N
|v_j(\gamma')|^2\\
&= \#(\Gamma) \sum_{j=1}^N
|v_j(\gamma')|^2.
\end{align*}
Since $\#(\Gamma) = \infty$, it follows that either
$N=0$ or $N=\infty$.
\end{proof}

\begin{cor}
Let $D$ be as in  Lemma 4.1 above.
Then the {\em essential spectrum}
of $D$ coincides with the spectrum of $D$, and so it is unbounded.
\end{cor}

\begin{proof}
By the Lemma above, we conclude that $spec_{ess}(D)$
and $spec(D)$ coincide. Since $D$ is an unbounded self-adjoint
operator, it is a standard fact that $spec(D)$ is unbounded
cf. \cite{Gi}, yielding the result.
\end{proof}

Note that in general the spectral projections of $D$,
$E_\lambda \not\in {C^*} (\mathcal G, \sigma)$. However
one has

\begin{prop}[Sunada, Bruning-Sunada]
Let $D$ be as in  Lemma 4.1 above.
If $\lambda_0 \not\in {\text{spec}}(D)$,
then $E_{\lambda_0} \in {C^*} (\mathcal G, \sigma)$.
\end{prop}

\begin{proof}
Firstly, there is a standard reduction to the case when $D$ is positive
and of even order $d\ge 2$ cf. \cite{BrSu}, so we will assume this without
loss of generality.
By a result of Greiner, see also Bruning-Sunada
\cite{BrSu}, there are off diagonal estimates for the Schwartz kernel of the
heat operator $e^{-t D}$
$$
|k_t(x,y)| \le C_1 t^{-n/d}exp\left(-C_2 d(x,y)^{d/(d-1)}t^{-1/(d-1)}\right)
$$
for some positive contants $C_1, C_2$ and for $t>0$ in any compact interval.
Since the volume growth of a orbifold covering space is at most exponential,
we see in particular that
$|k_t(x,y)|$ is $L^1$ in both variable seperately, so that
$$
e^{-t D} \in {C^*} (\mathcal G, \sigma).$$Note that
$\chi_{[0,e^{-t\lambda}]}(D) = \chi_{[0,\lambda]}(e^{-t D})$ Let $t=1$ and
$\lambda_1 = -\log\lambda_0$. Then $\lambda_1 \not\in
{\text{spec}}(e^{- D})$ and
$$
\chi_{[0,\lambda_1]}(e^{- D}) = \phi (e^{- D})
$$
where $\phi$ is a compactly supported smooth function, $\phi \cong 1$
on $[0,\lambda_1]$ and $\phi \cong 0$ on the remainder of the spectrum.
Since ${C^*} (\mathcal G, \sigma)$ is closed under the continuous
functional calculus, it follows that $\phi (e^{- D}) \in
{C^*} (\mathcal G, \sigma)$, that is $E_{\lambda_0}  \in
{C^*} (\mathcal G, \sigma)$.
\end{proof}

We will now recall the definition of the Kadison constant of a twisted group
$C^*$-algebra.
The {\em Kadison constant} of ${C}^*_r(\Gamma, \sigma)$ is defined by:
$$
C_\sigma(\Gamma) = \inf\{ \tr(P) : P \ \ {\mbox{is a
non-zero projection in}} \ \
{C}^*_r(\Gamma, \sigma) \otimes \mathcal K\}.
$$

\begin{prop} Let $\Gamma$ be a Fuchsian group of signature
$(g, \nu_1, \ldots, \nu_n)$.
Let $\sigma$ be a multiplier on
$\Gamma$ such that $\delta(\sigma) = 0$,  and
$2\pi\theta = \langle[\sigma], [\Gamma]\rangle \in (0,1]$
be the result of pairing the cohomology class of
$\sigma$ with the fundamental class of $\Gamma$. If $\theta$ is
rational, then
the spectrum of any $(\Gamma, \bar\sigma)$-invariant
elliptic differential operator $D$ on $\M$ has a band
structure, in the
sense that the intersection of the resolvent set with any compact interval
in $\mathbb R$
has only a finite number of components. Here $\Gamma\to \M\to M$
is the universal orbifold covering of a compact good orbifold $M$ with
orbifold fundamental group $\Gamma$.
In particular, the intersection of
$spec(D)$ with any compact interval in $\mathbb R$ is never a
Cantor set.
\end{prop}

\begin{proof}
By Proposition 4.1 and Theorem 3.3, it follows that one has the estimate
$C_\sigma(\Gamma)
\ge 1/q >0$. Then one applies the main result in Br\"uning-Sunada
\cite{BrSu} to
deduce the proposition.
\end{proof}

In words, we have shown that
whenever the multiplier is rational, then
the spectrum of a projectively periodic elliptic
operator is the union of countably many (possibly degenerate)
closed intervals, which can only accumulate at infinity.

Recall the important $\Gamma$-invariant
elliptic differential operator, which is called the
Schr\"odinger operator
$$
H_{V} = \Delta +V
$$
where $\Delta$ denotes the Laplacian on functions on
$\M$ and $V$ is a $\Gamma$-invariant function on $\M$.
It is known that the Baum-Connes conjecture is true for
all amenable discrete subgroups of a connected Lie group
and also for discrete subgroups of $SO(n,1)$ \cite{Kas1}
and $SU(n,1)$ \cite{JuKas}. For all these groups $\Gamma$,
it follows that the Kadison constant $C_1(\Gamma)$ is positive.
Therefore we see by the arguments above that the
spectrum of the periodic elliptic
operator $H_{V}$ is the union of countably many (possibly degenerate)
closed intervals, which can only accumulate at infinity.
This gives evidence for the following

\begin{conj*}[The Generalized Bethe-Sommerfeld conjecture]
The spectrum of any $\Gamma$-invariant Schr\"odinger operator
$H_{V}$ has
only a {\em finite} number of bands, in the
sense that the intersection of the resolvent set with $\mathbb R$
has only a finite number of components.
\end{conj*}

We remark that the Bethe-Sommerfeld conjecture has been
proved completely by Skriganov \cite{Skri} in the Euclidean case.

This leaves open the question of whether there are
$(\Gamma, \bar\sigma)$-invariant
elliptic differential operators $D$ on $\J$
with Cantor spectrum
when $\theta$ is irrational. In the Euclidean case,
this is usually known
as the {\em Ten Martini Problem}, and is to date,
not completely solved,
though much progress has been made (cf. \cite{Sh}). We pose a
generalization
of this problem to the hyperbolic case (which also includes the
Euclidean case):

\begin{conj*}[The Generalized Ten Dry Martini Problem]
Let $\sigma$ be a multiplier on
$\Gamma$ such that $\delta(\sigma) = 0$,  and
$2\pi\theta = \langle[\sigma], [\Gamma]\rangle \in (0,1]$
be the result of pairing the cohomology class of
$\sigma$ with the fundamental class of $\Gamma$. If $\theta$ is
irrational, then
there is a $(\Gamma, \bar\sigma)$-invariant
elliptic differential operator $D$ on $\J$ which has
a Cantor set type spectrum, in the sense that the intersection of
$spec(D)$ with some compact interval in
$\mathbb R$ is a Cantor set.
\end{conj*}

\subsection{On the classification of twisted group $C^*$-algebras}
We will now use the range of the trace Theorem 3.3, to give a complete
classification,
up to isomorphism, of the twisted group $C^*$-algebras ${C}^*(\Gamma,
\sigma)$, where we assume as before that $\delta(\sigma) = 0$.

\begin{prop}[Isomorphism classification of twisted group $C^*$-algebras]
Let $\sigma, \sigma' \in $ \\ $H^2(\Gamma, \mathbb R/ \mathbb Z)$ be
multipliers on
$\Gamma$ satisfying $\delta(\sigma) = 0 = \delta(\sigma')$, and
$2\pi\theta = <\sigma, [\Gamma]> \in (0,1]$,
$2\pi\theta' = <\sigma', [\Gamma]> \in (0,1]$
be the result of pairing
$\sigma$, $\sigma'$ with the fundamental class of $\Gamma$.
Then ${C}^*(\Gamma, \sigma) \cong {C}^*(\Gamma, \sigma')$ if and only if
$\theta' \in \{(\theta + \sum_{i=1}^n \beta_i/\nu_i) \quad \text{mod} 1,
(1-\theta + \sum_{i=1}^n \beta_i/\nu_i) \quad \text{mod} 1\}$, where
$0\le \beta_i \le \nu_i-1 \quad \forall i = 1,\ldots ,n$.
\end{prop}

\begin{proof}
 Let $\tr$ and $\tr'$ denote the canonical traces on ${C}^*(\Gamma, \sigma)$
and ${C}^*(\Gamma, \sigma')$ respectively. Let
$$
\phi : {C}^*(\Gamma, \sigma) \cong {C}^*(\Gamma, \sigma')
$$
be an isomorphism, and let
$$
\phi_* : K_0({C}^*(\Gamma, \sigma)) \cong K_0({C}^*(\Gamma, \sigma'))
$$
denote the induced map on $K_0$.
By Theorem 3.3, the range of the trace map on $K_0$ is
$$
[\tr] (K_0 ({C}^*(\Gamma, \sigma)) ) = \mathbb Z \theta + \mathbb Z +
\sum_{i=1}^n
\mathbb Z (1/\nu_i)
$$
and
$$
[\tr'] (K_0 ({C}^*(\Gamma, \sigma')) ) = \mathbb Z \theta' + \mathbb Z
+\sum_{i=1}^n
\mathbb Z (1/\nu_i).
$$
Therefore if $\theta$ is irrational, then $\mathbb Z \theta + \mathbb Z +
\sum_{i=1}^n
\mathbb Z (1/\nu_i) = \mathbb Z \theta' + \mathbb Z + \sum_{i=1}^n
\mathbb Z (1/\nu_i) $ implies that $\theta'$ is also irrational and
that $\theta \pm \theta' \in \mathbb Z + \sum_{i=1}^n
\mathbb Z (1/\nu_i) $.
Since $\theta, \theta' \in (0,1]$, one deduces that $\theta' \in \{(\theta
+ \sum_{i=1}^n \beta_i/\nu_i) \text{mod} 1,
(1-\theta + \sum_{i=1}^n \beta_i/\nu_i) \text{mod} 1\}$, where
$0\le \beta_i \le \nu_i-1 \quad \forall i = 1,\ldots ,n$. Virtually the
same argument holds when $\theta$ is rational, but one argues in $K$-theory
first, and applies the trace only at the final step.

First observe that a diffeomorphism $C: \Sigma_{g'} \to \Sigma_{g'}$
lifts to a diffeomorphism $C'$ of $\mathbb H$ such that $C' \Gamma {C'}^{-1}
= \Gamma$, i.e. it defines an automorphism of $\Gamma$. Recall that the finite
group $$
G = \left\{ C_i : \quad C_i^{\nu_i} = 1 \quad \forall i = 1, \ldots, n
\right\}
$$
acts on $\Sigma_{g'}$ with quotient $\Sigma(g, \nu_1, \ldots, \nu_n)$. By
the observation above, we see that $G$ also acts as automorphisms of
$\Gamma$. Now $C_i[\Gamma] = \lambda_i[\Gamma]$, where
$\lambda_i \in \mathbb C$. Since
$C_iC_j[\Gamma] = \lambda_i\lambda_j[\Gamma]$ and $C_i^{\nu_i} = 1$, it
follows that $\lambda_i$ is an $\nu_i^{th}$ root of unity, i.e.
$\lambda_i = e^{2\pi \sqrt{-1} (1/\nu_i)}$.
Let $C \in G$, i.e. $C = \prod_{i=1}^n C_i^{\beta_i}$.
We evaluate
$<C^*[\sigma], [\Gamma]> = <[\sigma], C_*[\Gamma]>
={<[\prod_{i=1}^n\lambda_i^{\beta_i}\sigma], [\Gamma]>}
= \theta + \sum_{i=1}^n \beta_i/\nu_i$. As in subsection 3.1 we see that
$C^*[\sigma] =
[\prod_{i=1}^n\lambda_i^{\beta_i}\sigma]\in
\ker\delta \subset H^2(\Gamma, U(1))$.
Therefore the automorphism $C_*$ of $\Gamma$ induces an isomorphism of
twisted
group $C^*$-algebras
$$
{C}^*(\Gamma, \sigma) \cong {C}^*(\Gamma, C^*\sigma)\cong
{C}^*(\Gamma, \lambda\sigma).
$$
where $\lambda = \prod_{i=1}^n\lambda_i^{\beta_i}$.

Now let $\psi : \Sigma_{g'} \to \Sigma_{g'}$ be an orientation reversing
diffeomorphism.
Then as observed earlier,
$\psi$ induces an automorphism $\psi_*:\Gamma \to \Gamma$ of $\Gamma$.
We evaluate
$<\psi^*[\sigma], [\Gamma]> = <[\sigma], \psi_*[\Gamma]>
={\overline{<[\sigma], [\Gamma]>}}=
<[\bar\sigma], [\Gamma]>$,
since $\psi$ is orientation
reversing. As in subsection 3.1 we see that $\psi^*[\sigma] =
[\bar\sigma]\in \ker\delta \subset H^2(\Gamma, U(1))$.
Therefore the automorphism $\psi_*$ of $\Gamma$ induces an isomorphism of
twisted
group $C^*$-algebras
$$
{C}^*(\Gamma, \sigma) \cong {C}^*(\Gamma, \psi^*\sigma)\cong
{C}^*(\Gamma, \bar\sigma).
$$

Therefore if
$\theta' \in \{(\theta + \sum_{i=1}^n \beta_i/\nu_i) \quad \text{mod} 1,
(1-\theta + \sum_{i=1}^n \beta_i/\nu_i) \quad \text{mod} 1\}$, where
$0\le \beta_i \le \nu_i-1 \quad \forall i = 1,\ldots ,n$,
one has
${C}^*(\Gamma, \sigma) \cong {C}^*(\Gamma, \sigma')$, completing the
proof of
the proposition.
\end{proof}

\section{Range of the higher trace on $K$-theory}

In this section, we compute the range of the
2-trace $\tr_c$ on $K$-theory of the twisted group $C^*$ algebra,
where $c$ is a 2-cocycle on the group,
generalising the work of \cite{CHMM}.
Suppose as before that $\Gamma$ is a discrete, cocompact subgroup of
$PSL(2,\R)$
of signature $(g, \nu_1, \ldots, \nu_n)$.
That is, $\Gamma$ is the orbifold fundamental group of a
compact hyperbolic orbifold $\Sigma(g,\nu_1,\ldots,\nu_n)$
of signature $(g,\nu_1,\ldots,\nu_n)$.  Then for any multiplier
$\sigma$ on $\Gamma$ such that $\delta(\sigma) =0$, one has the
\emph{twisted Kasparov isomorphism},
\[
   \mu_\sigma : K^\bullet_{orb} (\Sigma(g,\nu_1,\ldots,\nu_n)) \to
K_\bullet (C^*_r(\Gamma,\sigma))
\]
is defined as in the previous
section.
We note in the following section  (using a
result of \cite{Ji})
that given any
projection $P$ in $C^*_r(\Gamma,\sigma)$ there is
both a projection $\tilde P$ in the same $K_0$ class but lying in
a dense subalgebra, stable under the holomorphic functional calculus,
and a Fredholm module for this dense subalgebra, which may be
be paired with $\tilde P$
to obtain an analytic index. On the other hand, by the results of the current
section,
given any such projection $P$ there
is a topological index that we can associate to it. The main result we prove
here is that the range of the 2-trace $\tr_c$ on $K$-theory of the
twisted group $C^*$ algebra is always an integer multiple
of a rational number. This will enable us to compute the range of values
of the Hall conductance in the quantum Hall effect on hyperbolic space,
generalising the results in \cite{CHMM}.

The first step in the proof is to show that
given an additive
 group cocycle $c\in Z^2(\Gamma)$ we may define canonical pairings
with $K^0(\Sigma(g,\nu_1,\ldots,\nu_n))$ and $K_0(C^*_r(\Gamma, \sigma))$
which
are related by the twisted Kasparov isomorphism, by
generalizing some of the results of Connes and
Connes-Moscovici to the twisted case.
The group 2-cocycle $c$ may be
regarded as a skew symmetrised function on
$\Gamma\times\Gamma\times\Gamma$, so that
 we can use the results in section 2 to obtain a cyclic
2-cocycle $\tr_c$ on $\mathbb{C}(\Gamma, \sigma)\otimes \mathcal{R}$
by defining:
\[
   \tr_c(f^0\otimes r^0, f^1\otimes r^1, f^2\otimes r^2) =
{\mbox{Tr}}(r^0r^1r^2)
      \sum_{g_0g_1g_2=1} f^0(g_0) f^1(g_1) f^2(g_2) c(1, g_1, g_1g_2)
\sigma(g_1, g_2).
\]
Note that
 $\tr_c$ extends to $\mathbb{C}(\Gamma, \sigma) \otimes \mathcal{L}^2$, (where
$\mathcal{L}^2$ denotes Hilbert-Schmidt operators) and
since $\mathbb{C}(\Gamma, \sigma) \otimes
\mathcal{R} \subset \mathbb{C}(\Gamma, \sigma) \otimes \mathcal{L}^2$   by
the pairing
theory of \cite{Co} one gets an additive map
\[
   [\tr_c] : K_0(\mathcal{R}(\Gamma, \sigma))\to\mathbb{R}.
\]
Explicitly, $[\tr_c]([e]-[f]) = \widetilde\tr_c(e,\cdots,e) -
\widetilde\tr_c(f,\cdots,f)$, where $e,f$ are idempotent matrices with
entries in $(\mathcal{R}(\Gamma, \sigma))^\sim = $ unital
algebra obtained by adding the identity to
$\mathcal{R}(\Gamma, \sigma)$ and $\widetilde\tr_c$
denotes the
canonical extension of $\tr_c$ to $(\mathcal{R}(\Gamma, \sigma))^\sim$. Let \
$\widetilde{\npartial_{\mathcal{E}}^+}\otimes \nabla$ \  be
the Dirac operator defined in the previous section, which is
invariant under the projective action of the fundamental group
defined by $\sigma$.
Recall that by definition, the $(c,\Gamma, \sigma)$-index of
\ $\widetilde{\npartial_{\mathcal{E}}^+}\otimes \nabla$\ is
\[
   \Ind_{c, \Gamma, \sigma}(\widetilde{\npartial_{\mathcal{E}}^+}\otimes
\nabla)=
 [\tr_c](\Ind_{(\Gamma, \sigma)} (\widetilde{\npartial_{\mathcal{E}}^+}
\otimes\nabla)) = \langle[\tr_c], \mu_\sigma([\mathcal E])\rangle
\in\mathbb{R}.
\]
It only depends on the cohomology class $[c]\in H^2(\Gamma)$, and it is
linear with respect to $[c]$.  We assemble this to give the following theorem.

\begin{thm}
Given $[c] \in H^2(\Gamma)$ and $\sigma \in H^2(\Gamma, U(1))$ a
multiplier on $\Gamma$, there is a canonical additive map
\[
   \langle [c],\ \ \rangle :
K^0_{orb}(\Sigma(g,\nu_1,\ldots,\nu_n))\to\mathbb{R},
\]
which is defined as
\[
   \langle [c], [\mathcal{E}]\rangle =
\Ind_{c, \Gamma, \sigma}(\widetilde{\npartial_{\mathcal{E}}^+}\otimes \nabla)=
 [\tr_c](\Ind_{(\Gamma, \sigma)} (\widetilde{\npartial_{\mathcal{E}}^+}
\otimes\nabla)) = \langle[\tr_c], \mu_\sigma([\mathcal E])\rangle
\in\mathbb{R}.
\]
Moreover, it is linear with respect to $[c]$.
\end{thm}

The {\em area cocycle} $c$ of the Fuchsian group $\Gamma$ is a canonically
defined 2-cocycle
on $\Gamma$ that is defined as follows.
Firstly, recall that
there is a well known area 2-cocycle on $PSL(2, \R)$cf. \cite{Co2}
defined as follows:
$PSL(2, \R)$ acts on $\mathbb H$ such that $\mathbb H \cong PSL(2, \R)/SO(2)$.
Then $c(g_1, g_2) = \text{Area}(\Delta(o, g_1.o, {g_2}^{-1}.o)) \in \R$,
where $o$
denotes an origin in $\mathbb H$ and
$\text{Area}(\Delta(a,b,c))$ denotes the hyperbolic area of the
geodesic triangle in $\mathbb H$
with vertices at $a, b, c \in \mathbb H$. Then the restriction of
$c$ to the subgroup $\Gamma$ is the  area cocycle $c$ of $\Gamma$.

\begin{cor}
Let $c,\ [c]\in H^2(\Gamma)$, be the area cocycle.  Then one has
\[
   \langle [c], [{\mathcal{E}}]\rangle = \phi\rank\mathcal{E}
      \in\phi\mathbb{Z}.
\]
where $ \phi = 2(g-1) + (n-\nu)\in {\mathbb Q}$ and
$\quad \nu = \sum_{j=1}^n 1/\nu_j$.
\end{cor}

\begin{proof}
By our generalization of the Connes-Moscovici higher index theorem
\cite{CM} to the twisted case of
elliptic operators on an orbifold covering space and that are invariant
under the projective action of the orbifold
fundamental group defined by $\sigma$,
cf. section 2,
 one has
$$
   [\tr_c](\Ind_{(\Gamma, \sigma)}
(\widetilde{\npartial_{\mathcal{E}}^+}\otimes
\nabla)) =
       \frac{1}{2\pi\#(G)}\int_{\Sigma_{g'}} \hat{A}(\Omega)
\tr(e^{R^{\mathcal{E}}}) e^{\omega}
\psi^*(\tilde c),\eqno{(5.1)}
$$
where $\psi: \Sigma_{g'} \to \Sigma_{g'}$ is the lift of the
map $f:\Sigma(g, \nu_1, \ldots \nu_n) \to \Sigma(g, \nu_1, \ldots \nu_n)$
(since $\underline B\Gamma = \Sigma(g, \nu_1, \ldots \nu_n)$ in this case)
which is the classifying map of the orbifold universal
cover (and which in this case is the identity map) and $[\tilde c]$
degree 2 cohomology class on $\Sigma_{g'}$ that
is the lift of $c$ to $\Sigma_{g'}$.
We next simplify the right hand
side of $(5.1)$ using the fact that $\hat{A}(\Omega) = 1$ and that
\begin{align*}
   \tr(e^{R^{\mathcal{E}}}) &= \rank \mathcal{E} + \tr({R^{\mathcal{E}}}), \\
   \psi^*(\tilde c) &= \tilde c,\\
   e^{\omega} &= 1 + {\omega} .
\end{align*}
We obtain
\[
   [\tr_c](\Ind_{(\Gamma, \sigma)}
(\widetilde{\npartial_{\mathcal{E}}^+}\otimes
\nabla)) =
      \frac{\rank \mathcal E}{2\pi\#(G)}
\langle [\tilde c], [\Sigma_{g'}] \rangle .
\]
When $c,\ [c]\in H^2(\Gamma)$, is the area 2-cocycle, then
$\tilde c$ is merely the restriction of the area cocycle on $PSL(2, \R)$
to the subgroup $\Gamma_{g'}$. Then one has
\[
   \langle [\tilde c], [\Sigma_{g'}] \rangle
= -2\pi\chi(\Sigma_{g'}) = 4\pi(g'-1).
\]
The corollary now follows from Theorem 5.1 above and the fact that
$\quad g' = 1 + \frac{\#(G)}{2}(2(g-1) + (n-\nu))$, $\quad \nu =
\sum_{j=1}^n 1/\nu_j$
and $\#(G) = 1+\sum_{j=1}^n (\nu_j -1)$.
\end{proof}

We next describe the canonical pairing of $K_0(C^*_r(\Gamma, \sigma))$,
given $[c]\in H^2(\Gamma)$. Since $\Sigma(g,\nu_1,\ldots,\nu_n)$ is
negatively curved, we
know from \cite{Ji} that
\[
   {\mathcal R}(\Gamma, \sigma) = \left\{ f : \Gamma \to \mathbb{C} \mid
      \sum_{\gamma\in\Gamma} |f(\gamma)|^2 (1+l(\gamma))^k < \infty
      \mbox{ for all } k\ge 0\right\},
\]
where $l:\Gamma \to \mathbb{R}^+$ denotes the length function, is a dense
and spectral invariant subalgebra of $C_r^*(\Gamma, \sigma)$.  In
particular it is closed under the smooth functional calculus, and is
known as the algebra of rapidly decreasing $L^2$ functions on $\Gamma$.
By a theorem of \cite{Bost}, the inclusion map ${\mathcal R}(\Gamma,
\sigma)\subset
C^*_r(\Gamma, \sigma)$ induces an isomorphism
\[
   K_j({\mathcal R}(\Gamma, \sigma)) \cong K_j(C_r^*(\Gamma, \sigma)),\quad
j=0,1.
\]
As  $\Sigma(g,\nu_1,\ldots,\nu_n) = \underline B\Gamma$
is a negatively curved orbifold,
we know (by \cite{Mos} and \cite{Gr}) that degree 2 cohomology classes in
$H^2(\Gamma)$ have
\emph{bounded} representatives i.e.\  bounded 2-cocycles on $\Gamma$.
Let $c$ be a bounded 2-cocycle on $\Gamma$.  Then it defines a cyclic
2-cocycle $\tr_c$ on the twisted group algebra
${\mathbb C}(\Gamma, \sigma)$, by a slight modification of
the standard formula \cite{CM}, (\cite{Ma} for the general case)
\[
   \tr_c(f^0, f^1, f^2) = \sum_{g_0g_1g_2=1} f^0(g_0) f^1(g_1) f^2(g_2)
      c(1, g_1, g_1g_2)\sigma(g_1, g_2).
\]
Here $c$ is assumed to be skew-symmetrized. Since the only difference
with the expression obtained in \cite{CM} is $\sigma(g_1, g_2)$, and since
$|\sigma(g_1, g_2)| = 1$, we can use Lemma 6.4, part (ii) in \cite{CM}
and the assumption that $c$ is bounded, to obtain the necessary
estimates which show that in fact $\tr_c$ extends continuously to the bigger
algebra ${\mathcal R}(\Gamma, \sigma)$. This induces an additive map in
$K$-theory as before:
\begin{gather*}
   [\tr_c] : K_0({\mathcal R}(\Gamma, \sigma))\to\mathbb{R} \\
   [\tr_c]([e] - [f]) = \widetilde\tr_c(e,\cdots, e)
      -\widetilde\tr_c(f,\cdots,f),
\end{gather*}
where $e,f$ are idempotent matrices with entries in $({\mathcal R}(\Gamma,
\sigma))^\sim$ (the
unital algebra associated to ${\mathcal R}(\Gamma, \sigma)$) and
$\widetilde\tr_c$ is
the canonical extension of $\tr_c$ to $({\mathcal R}(\Gamma, \sigma))^\sim$.
Observe that the twisted Kasparov map is merely
$$
\mu_\sigma([\mathcal{E}]) = j_*(\Ind_{(\Gamma, \sigma)}
(\widetilde{\npartial^+_{\mathcal{E}}}\otimes \nabla) ) \in
K_0({C}^*(\Gamma, \sigma)).
$$
Here $j: {\mathbb C}(\Gamma, \sigma)\otimes \mathcal R \to
{C}^*(\Gamma, \sigma)\otimes \mathcal K$
is the natural inclusion map, and $j_* : K_0 ( {\mathbb C}(\Gamma,
\sigma)\otimes \mathcal R) \to K_0 ({C}^*(\Gamma, \sigma))$ is the
induced map in $K$-theory. Therefore one has the equality
\[
   \langle [c], \mu_\sigma^{-1}[P]\rangle = \langle [\tr_c], [P] \rangle
\]
for any $[P]\in K_0({\mathcal R}(\Gamma, \sigma)) \cong K_0(C_r^* (\Gamma,
\sigma))$.
Using the previous corollary, one has

\begin{thm}[Range of the higher trace on $K$-theory]
Let $c$ be the area 2-cocycle on $\Gamma$.  Then $c$ is known to
be a bounded 2-cocycle, and one has
\[
   \langle [\tr_c], [P] \rangle = \phi(\rank\mathcal{E}^0 -
      \rank\mathcal{E}^1)\in\phi\mathbb{Z},
\]
where $ \phi = 2(g-1) + (n-\nu)\in {\mathbb Q}$
and $\quad \nu = \sum_{j=1}^n 1/\nu_j$.
Here $[P]\in K_0({\mathcal R}(\Gamma, \sigma)) \cong K_0(C_r^*(\Gamma,
\sigma))$,
and $\mathcal{E}^0, \ {\mathcal E}^1$ are orbifold vector bundles over
$\Sigma(g, \nu_1, \ldots,\nu_n)$ such that
\[
   \mu_\sigma^{-1}([P]) = [\mathcal{E}^0] -
      [{\mathcal E}^1] \in K^0_{orb}(\Sigma(g,\nu_1,\ldots,\nu_n)).
\]
In particular, the range of the the higher trace on $K$-theory is
$$
[\tr_c]( K_0 (C^*(\Gamma,\sigma))) = \phi\Z.
$$
\end{thm}

Note that
$\phi$ is in general only a {\em rational} number and we will give examples
to show that this is the case; however it is
an {\em integer} whenever the orbifold is smooth, i.e. whenever $1=\nu_1=\ldots
=\nu_n$,which was the case that was considered in \cite{CHMM}.
We will apply this result in section 6
to compute the range of values the Hall conductance in the quantum Hall effect
on the hyperbolic plane,
for orbifold fundamental groups, extending the results in \cite{CHMM}.

\subsection{Examples of orbifolds with fractional $\phi$} In \cite{Bro},
Broughton has listed all the good two dimensional orbifolds
which are quotients of Riemann surfaces $\Sigma_{g'}$ with genus
$g' = 2\, \text{or} \, 3$. Using his explicit classification, we will give
several
examples showing that the number $\phi$, as in Theorem 5.3,
can be a fraction. Even when $g'=2$, there are several examples
from Table 4, \cite{Bro}, but we will focus on one of these, viz.
2.k.2, where a dihedral group of order 6 acts on the genus two
Riemann surface giving rise to an orbifold of signature
$(g, \nu_1, \nu_2, \nu_3, \nu_4) = (0, 2,2,3,3)$. It follows that
$\phi = 1/3$ in this case, and that the range of values
the Hall conductance in the quantum Hall effect
on the hyperbolic plane, for this particular orbifold fundamental
group is $\Z(1/3)$. One can list all the possible
denominators that can occur for $\phi$ by looking through
the Tables 4 and 5 in \cite{Bro}, which will be used again
to determine all the range of values
the Hall conductance in the quantum Hall effect
on the hyperbolic plane, and will be studied further in
section 6.

\section{The Area cocycle, the hyperbolic Connes-Kubo formula and
the Quantum Hall Effect}

In this section we prove a generalisation of the results in
\cite{CHMM} on the Quantum Hall Effect on the hyperbolic plane, where
we now allow the discrete group to have torsion. We will only discuss
the discrete model, as the discussion for the continuous model is similar,
as shown in \cite{CHMM}. We will first derive the discrete analogue
of the hyperbolic Connes-Kubo formula for the Hall conductance 2-cocycle,
which was derived in the continuous case in \cite{CHMM}. We then relate it
to the Area 2-cocycle on the twisted group algebra of the
discrete Fuchsian group, and we show that these define the same
cyclic cohomology class.
This enables us to use the results of the previous section to
show that the Hall conductance
has plateaus at all energy levels belonging to
any gap in the spectrum of the Hamiltonian, where it
is now shown shown to be equal to an integral multiple of
a {\em fractional} valued
invariant.  Moreover the set of possible
denominators is finite
and has been explicitly determined. It is plausible that this
might shed light on the mathematical mechanism
responsible for fractional quantum numbers.

The graph that we consider is the Cayley graph
of the Fuchsian group $\Gamma$ of signature
$(g, \nu_1, \ldots, \nu_n)$, which acts freely on
the complement of a countable set of points in the
hyperbolic plane.
The Cayley graph embeds in the hyperbolic plane as follows.
Fix a base point $u \in \mathbb H$ such that the stabilizer (or isotropy
subgroup)
at $u$ is trivial and consider the orbit of the $\Gamma$ action through
$u$. This gives the vertices of the graph. The
edges of the graph are geodesics constructed as follows.
Each element of the group $\Gamma$
may be written as a word of minimal length in the
generators of $\Gamma$ and their inverses. Each generator
and its inverse determines a unique geodesic emanating from a vertex $x$
and these form the edges  of the
graph. Thus each word $x$ in the generators determines a
piecewise geodesic path from $u$ to $x$. Note that if
$\Gamma$ has elliptic elements, then the Cayley graph of
$\Gamma$ has (geodesic) loops i.e. it is not a tree.

Recall that the {\em area cocycle} $c$ of the Fuchsian group
$\Gamma$ is a canonically defined 2-cocycle
on $\Gamma$ that is defined as follows.
Firstly, recall that
there is a well known area 2-cocycle on $PSL(2, \R)$, cf. \cite{Co2},
defined as follows:
$PSL(2, \R)$ acts on $\mathbb H$ such that
$\mathbb H \cong PSL(2, \R)/SO(2)$.
Then $c(\gamma_1, \gamma_2) = \text{Area}(\Delta(o, \gamma_1.o,
{\gamma_2}^{-1}.o)) \in \R$,
where $o$
denotes an origin in $\mathbb H$ and
$\text{Area}(\Delta(a,b,c))$ denotes the hyperbolic area of the
geodesic triangle in $\mathbb H$
with vertices at $a, b, c \in \mathbb H$. Then the restriction of
$c$ to the subgroup $\Gamma$ is the  area cocycle $c$ of $\Gamma$.

This area cocycle defines in a canonical way a cyclic 2-cocycle
$\tr_c$ on the group algebra $\C(\Gamma,\sigma)$ as follows;
$$
\tr_c(a_0,a_1,a_2) =
\sum_{\gamma_0\gamma_1\gamma_2 =1}
a_0(\gamma_0) a_1 (\gamma_1) a_2 (\gamma_2) c(\gamma_1, \gamma_2)
\sigma(\gamma_1, \gamma_2)
$$

We will now describe the hyperbolic Connes-Kubo formula for the Hall
conductance in the Quantum Hall Effect. Let $\Omega_j$ denote the (diagonal)
operator on $\ell^2(\Gamma)$ defined by
$$
\Omega_j f(\gamma) = \Omega_j(\gamma)f(\gamma) \quad \forall f\in
\ell^2(\Gamma)
\quad \forall \gamma\in \Gamma
$$
where
$$
\Omega_j(\gamma) = \int_o^{\gamma.o} a_j \quad j=1,\ldots ,2g
$$
and where $\, \{a_j\} \quad j=1,\ldots ,2g$  is the lift to $\mathbb H$ of
a symplectic basis of harmonic 1-forms on the Riemann surface of genus $g$
underlying the orbifold
$\Sigma(g, \nu_1,\ldots,\nu_n)$.

For $j=1,\ldots ,2g$, define the derivations $\delta_j$ on
$\mathcal{R}(\Gamma,\sigma)$
as being the commutators $\delta_j a = [\Omega_j, a]$. A simple calculation
shows that
$$
\delta_j a (\gamma) = \Omega_j(\gamma)
a(\gamma) \quad \forall a\in \mathcal{R}(\Gamma,\sigma)
\quad \forall \gamma\in \Gamma.
$$
Note that these are not inner derivations, and also that we have the simple
estimate
$$
|\Omega_j(\gamma)|\le ||a_j||_{(\infty)} d(\gamma.o, o)
$$
where $d(\gamma.o, o)$ and the distance in the word metric on the group
$\Gamma$, $d_\Gamma (\gamma, 1)$ are equivalent. This then yields the
estimate
$$
|\delta_j a (\gamma)| \le C_N d_\Gamma (\gamma, 1)^{-N}\quad \forall N\in
\mathbb N
$$
i.e  $\delta_j a \in \mathcal{R}(\Gamma,\sigma) \quad \forall a
\in \mathcal{R}(\Gamma,\sigma)$. Note that since $\forall \gamma,
\gamma' \in \Gamma$, the difference
$\Omega_j(\gamma\gamma') - \Omega_j(\gamma')$ is a constant
independent of $\gamma'$, we see that $\Gamma$-equivariance is preserved.
For $j=1,\ldots ,2g$, define the cyclic 2-cocycles
$$
\tr^K_j (a_0,a_1,a_2) = \tr(a_0 (\delta_ja_1 \delta_{j+g} a_2
- \delta_{j+g} a_1 \delta_j a_2)).
$$
These are supposed to
give the Hall conductance for currents in the $(j+g)$th direction
which are induced by electric fields in the $j$th direction
cf. section 6, \cite{CHMM}.
Then the {\em hyperbolic Connes-Kubo formula} for the Hall
conductance is the cyclic 2-cocycle given by the sum
$$
\tr^K (a_0,a_1,a_2) =
\sum_{j=1}^g \tr^K_j (a_0,a_1,a_2).
$$

\begin{thm}[The Comparison Theorem]
$$
[\tr^K] =  [\tr_c] \in HC^2(\mathcal{R}(\Gamma,\sigma))
$$
\end{thm}

\noindent{\bf Sketch of Proof:}
Our aim is now to compare the two cyclic 2-cocycles
and to prove that they differ by a coboundary i.e.
$$
\tr^K (a_0,a_1,a_2) - \tr_c (a_0,a_1,a_2) =
b\lambda (a_0,a_1,a_2)
$$
for some cyclic 1-cochain $\lambda$ and where $b$ is the cyclic
coboundary operator. The key to this theorem is a
geometric interpretation
of the hyperbolic Connes-Kubo formula.

We begin with some calculations, to enable us to make this
comparison of the cyclic 2-cocycles.
$$
\tr^K (a_0,a_1,a_2) =
$$
$$
\sum_{j=1}^g \sum_{\gamma_0\gamma_1\gamma_2 =1}
a_0(\gamma_0) \left( \delta_j a_1(\gamma_1) \delta_{j+g} a_2(\gamma_2)
- \delta_{j+g} a_1(\gamma_1) \delta_j a_2(\gamma_2)\right)
\sigma(\gamma_0, \gamma_1) \sigma(\gamma_0\gamma_1, \gamma_2)
$$
$$
= \sum_{j=1}^g \sum_{\gamma_0\gamma_1\gamma_2 =1}
a_0(\gamma_0) a_1(\gamma_1) a_2(\gamma_2) \left(
\Omega_j(\gamma_1) \Omega_{j+g}(\gamma_2)
- \Omega_{j+g}(\gamma_1) \Omega_j(\gamma_2)
\right)\sigma(\gamma_1, \gamma_2)
$$
since by the cocycle identity for multipliers, one has
\begin{align*}
\sigma(\gamma_0, \gamma_1) \sigma(\gamma_0\gamma_1, \gamma_2)
& =  \sigma(\gamma_0,\gamma_1\gamma_2) \sigma(\gamma_1, \gamma_2)\\
& = \sigma(\gamma_0,\gamma_0^{-1}) \sigma(\gamma_1, \gamma_2) \quad
\text{since}
\quad \gamma_0\gamma_1\gamma_2 = 1\\
& = \sigma(\gamma_1, \gamma_2) \quad \text{since}
\quad\sigma(\gamma_0,\gamma_0^{-1}) =1.
\end{align*}
So we are now in a position to compare the two cyclic 2-cocycles.
Define $\Psi_j(\gamma_1, \gamma_2) =
\Omega_j(\gamma_1) \Omega_{j+g}(\gamma_2)
- \Omega_{j+g}(\gamma_1) \Omega_j(\gamma_2)$.
Let $\Xi: \mathbb H \to {\mathbb R}^{2g}$ denote the lift of the Abel-Jacobi
map. It is a symplectic map, since it is known to be a holomorphic embedding.
Therefore if $\omega$ and $\omega_J$ are their respective symplectic
2-forms, then one has $\Xi^*(\omega_J) = \omega$.
Then one has the following
geometric lemma

\begin{lemma}
$$
\sum_{j=1}^g  \Psi_j(\gamma_1, \gamma_2) = \int_{\Delta_E(\gamma_1, \gamma_2)}
\omega_J
$$
where $\Delta_E(\gamma_1, \gamma_2)$ denotes the Euclidean triangle with
vertices at $\Xi(o), \Xi(\gamma_1.o)$ and $ \Xi(\gamma_2.o)$, and
$\omega_J$ denotes the flat K\"ahler 2-form on the Jacobi variety. That is,
$\sum_{j=1}^g  \Psi_j(\gamma_1, \gamma_2)$ is equal to the
Euclidean
area of the Euclidean triangle $\Delta_E(\gamma_1, \gamma_2)$.
\end{lemma}

\begin{proof} We need to consider the expression
$$\sum_{j=1}^g  \Psi_j(\gamma_1, \gamma_2)
= \sum_{j=1}^g \Omega_j(\gamma_1) \Omega_{j+g}(\gamma_2)
- \Omega_{j+g}(\gamma_1) \Omega_j(\gamma_2).$$
Let $s$ denote the symplectic form on ${\mathbb R}^{2g}$ given by:
$$s(u,v)=\sum_{j=1}^g(u_jv_{j+g} -u_{j+g}v_j).$$
The so-called `symplectic area' of a triangle
with vertices $\Xi(o)=0,\Xi(\gamma_1.o),\Xi(\gamma_2.o)$ may be seen to be
$s(\Xi(\gamma_1.o),\Xi(\gamma_2.o))$. To appreciate this, however,
we need to utilise an argument from
\cite{GH}, pages 333-336.
In terms of the standard basis of ${\mathbb R}^{2g}$ (given in this case
by vertices in the integer period lattice arising from our choice of basis
of harmonic one forms) and corresponding coordinates
$u_1,u_2,\ldots u_{2g}$
the form $s$ is the two form on ${\mathbb R}^{2g}$ given by
$$\omega_J=\sum_{j=1}^g du_j\wedge du_{j+g}.$$
Now the  `symplectic area' of a triangle in ${\mathbb R}^{2g}$ with
 vertices $\Xi(o)=0,\Xi(\gamma_1.o),\Xi(\gamma_2.o)$
is given by integrating $\omega_J$ over the triangle and
a brief calculation reveals that this yields
$s(\Xi(\gamma_1.o),\Xi(\gamma_2.o))/2$, proving the lemma.
\end{proof}

We also observe that since $\omega = \Xi^* \omega_J$, one has
$$
c(\gamma_1, \gamma_2) = \int_{\Delta(\gamma_1, \gamma_2)}
\omega = \int_{\Xi(\Delta(\gamma_1, \gamma_2))}
\omega_J
$$
Therefore the difference
\begin{align*}
\sum_{j=1}^g  \Psi_j(\gamma_1, \gamma_2)  - c(\gamma_1, \gamma_2)
& = \int_{\Delta_E(\gamma_1, \gamma_2)}
\omega_J  - \int_{\Xi(\Delta(\gamma_1, \gamma_2))}
\omega_J\\
& = \int_{\partial \Delta_E(\gamma_1, \gamma_2)}
\Theta_J  - \int_{\partial\Xi(\Delta(\gamma_1, \gamma_2))}
\Theta_J
\end{align*}
where $\Theta_J$ is a 1-form on the universal cover of the
Jacobi variety such that $d\Theta_J = \omega_J$.  Therefore
one has
$$
\sum_{j=1}^g  \Psi_j(\gamma_1, \gamma_2)  - c(\gamma_1, \gamma_2)
= h(1, \gamma_1)  - h(\gamma_1^{-1}, \gamma_2)  +
h(\gamma_2^{-1}, 1)
$$
where $h(\gamma_1^{-1}, \gamma_2) =  \int_{\Xi(\ell (\gamma_1, \gamma_2))}
\Theta_J - \int_{m(\gamma_1, \gamma_2)}
\Theta_J $, where $\ell(\gamma_1, \gamma_2)$ denotes the
unique geodesic in $\mathbb H$ joining $\gamma_1.o$ and $\gamma_2.o$
and $m(\gamma_1, \gamma_2)$ is the straight line in the Jacobi variety
joining the points $\Xi(\gamma_1.o)$ and $\Xi(\gamma_2.o)$.
Since we can also write
$h(\gamma_1^{-1}, \gamma_2) = \int_{D(\gamma_1, \gamma_2)}
\omega_J
$, where $D(\gamma_1, \gamma_2)$ is a disk in the Jacobi variety
with boundary $\Xi(\ell(\gamma_1, \gamma_2))\cup m(\gamma_1, \gamma_2)$,
we see that $h$ is $\Gamma$-invariant.

We now define the cyclic 1-cochain $\lambda$ on $\mathcal{R}(\Gamma,\sigma)$
as
$$
\lambda(a_0, a_1) = \tr((a_{0})_h a_1) = \sum_{\gamma_0\gamma_1 =1}
h(1, \gamma_1) a_0(\gamma_0) a_1(\gamma_1)\sigma(\gamma_0, \sigma_1)
$$
where $(a_{0})_h$ is the operator on $\ell^2(\Gamma)$ whose matrix
in the canonical basis is
$h(\gamma_1, \gamma_2)a_0(\gamma_1\gamma_2^{-1})$.
Firstly, one has by definition
$$
b\lambda(a_0, a_1, a_2) =  \lambda(a_0 a_1, a_2)-
\lambda(a_0, a_1a_2)+\lambda(a_2a_0, a_1)
$$
We compute each of the terms seperately
\begin{align*}
\lambda(a_0 a_1, a_2) & =
\sum_{\gamma_0\gamma_1\gamma_2 =1}
h(1, \gamma_2) a_0(\gamma_0) a_1(\gamma_1)a_2(\gamma_2)
\sigma(\gamma_1, \gamma_2) \\
\lambda(a_0, a_1 a_2) & =
\sum_{\gamma_0\gamma_1\gamma_2 =1}
h(1, \gamma_1 \gamma_2) a_0(\gamma_0) a_1(\gamma_1)a_2(\gamma_2)
\sigma(\gamma_1, \gamma_2) \\
\lambda(a_2 a_0, a_1) & =
\sum_{\gamma_0\gamma_1\gamma_2 =1}
h(1, \gamma_1) a_0(\gamma_0) a_1(\gamma_1)a_2(\gamma_2)
\sigma(\gamma_1, \gamma_2)
\end{align*}
Now by $\Gamma$-equivariance, $h(1, \gamma_1\gamma_2)
= h(\gamma_1^{-1}, \gamma_2)$ and  $
h(1, \gamma_2) =  h(\gamma_2^{-1}, 1)$. Therefore one has
$$
b\lambda(a_0, a_1, a_2) =
$$
$$
\sum_{\gamma_0\gamma_1\gamma_2 =1}
 a_0(\gamma_0) a_1(\gamma_1)a_2(\gamma_2)
\left(h(\gamma_2^{-1}, 1) - h(\gamma_1^{-1}, \gamma_2)
+ h(1, \gamma_1)\right)
\sigma(\gamma_1, \gamma_2)
$$
Using the formula above, we see that
$$
b\lambda(a_0, a_1, a_2) =   \tr^K (a_0,a_1,a_2) - \tr_c (a_0,a_1,a_2).
$$

It follows from Connes pairing theory of cyclic cohomology and $K$-theory
\cite{Co2} and the Comparison Theorem above that

\begin{cor}
$$
\tr^K (P, P, P) = \tr_c (P, P, P)
$$
for all projections $P\in \mathcal{R}(\Gamma,\sigma)$.
\end{cor}

Recall that by the range of the higher trace Theorem 5.3, one has
$$
\tr_c (P, P, P) \in \phi\mathbb Z \eqno(6.1)
$$
for all projections $P\in \mathcal{R}(\Gamma,\sigma)$, where
$\phi = 2(g-1) + (n-\nu) \in \mathbb Q$.

Finally, suppose that we are given a very thin hyperboloid sample
of pure metal, with electrons situated along
the Cayley graph of $\Gamma$, and a very strong
magnetic field which is uniform and normal in direction to the sample.
Then at very low temperatures, close to absolute zero,
quantum mechanics dominates and the discrete model that is considered
here is a model of electrons moving on the
Cayley graph of $\Gamma$ which is
embedded in the hyperboloid.
The associated discrete Hamiltonian $H_\sigma$
for the electron in the magnetic field
is given by the Random Walk operator in the
projective $(\Gamma, \sigma)$ regular representation on the
Cayley graph of the group $\Gamma$. It is also known as the generalized
{\em Harper operator} and was first studied in this generalized
context in \cite{Sun} and also in \cite{CHMM}. We will see that
the Hamiltonian that we consider is in a natural way
the sum of
a free Hamiltonian and an interacting term.

The hyperbolic Connes-Kubo formula for the
Hall conductance $\sigma_E$ at the energy level $E$ is
defined as follows; let $P_E = \chi_{[0,E]}(H_\sigma)$ be the spectral
projections of the Hamiltonian to energy levels less than or
equal to $E$. Then if $E\not\in \text{spec}(H_\sigma)$, one can show that
$P_E \in \mathcal{R}(\Gamma,\sigma)$, and the
Hall conductance is defined as
$$
\sigma_E = \tr^K (P_E, P_E, P_E).
$$
As mentioned earlier, it measures the sum of the contributions
to the Hall conductance at the energy level $E$
for currents in the $(j+g)$th direction which are induced by
electric fields in the $j$th direction, cf. section 6 \cite{CHMM}.
When this is combined with equation $(6.1)$,
one sees that the Hall conductance takes on values
in $\phi \mathbb Z$  whenever
the energy level $E$ lies in a gap in the spectrum of the Hamiltonian
$H_\sigma$. In fact we notice that the Hall conductance is a constant function
of the energy level $E$ for all values of $E$ in the same gap
in the spectrum of the Hamiltonian. That is, the Hall conductance
has plateaus which are {\em integer} multiples of the fraction $\phi$
on the gap in the spectrum of the Hamiltonian.

We now give some details. Recall the left $\sigma$-regular representation
$$
(\lambda(\gamma) f)(\gamma') = f(\gamma^{-1}\gamma')
\sigma(\gamma', \gamma^{-1}\gamma')
$$
$\forall f\in \ell^2(\Gamma)$ and $\forall \gamma, \gamma' \in \Gamma$.
It has the property that
$$
\lambda(\gamma) \lambda(\gamma') = \sigma(\gamma, \gamma')
\lambda(\gamma\gamma')
$$
Let $S= \{A_j, B_j, A_j^{-1}, B_j^{-1}, C_i, C_i^{-1} : j= 1, \ldots, g, \quad
i=1, \ldots,n\}$ be a
symmetric set of generators for $\Gamma$. Then the Hamiltonian
is explicitly given as
\begin {align*}
H_\sigma & :\ell^2(\Gamma) \to \ell^2(\Gamma)\\
H_\sigma & = \sum_{\gamma\in S} \lambda(\gamma)
\end{align*}
and is clearly by definition a bounded self adjoint operator.
Notice that the Hamiltonian can be decomposed as a sum of a free
Hamiltonian containing the torsionfree generators
and an interaction term containing the torsion generators.
$$
H_\sigma  = H_\sigma^{free} + H_\sigma^{interaction}
$$
where
$$
H_\sigma^{free} = \sum_{\gamma\in S'} \lambda(\gamma) \qquad \text{and} \qquad
H_\sigma^{interaction}  =\sum_{\gamma\in S''} \lambda(\gamma)
$$
and where $S'= \{A_j, B_j, A_j^{-1}, B_j^{-1}: j= 1, \ldots, g\}$ and
$S''= \{ C_i, C_i^{-1} : i=1, \ldots,n\}$.

\begin{lemma}
If $E\not\in \text{spec}(H_\sigma)$, then $P_E \in \mathcal{R}(\Gamma,\sigma)$,
where $P_E = \chi_{[0,E]}(H_\sigma)$ is the spectral
projection of the Hamiltonian to energy levels less than or
equal to $E$.
\end{lemma}

\begin{proof}
Since $E\not\in \text{spec}(H_\sigma)$, then
$P_E = \chi_{[0,E]}(H_\sigma) = \phi(H_\sigma)$ for some smooth,
compactly supported function $\phi$. Now by definition, $H_\sigma \in
\C(\Gamma,\sigma) \subset  \mathcal{R}(\Gamma,\sigma)$,
and since $\mathcal{R}(\Gamma,\sigma)$ is closed under the smooth
functional calculus by the result of \cite{Ji}, it follows that
$P_E \in \mathcal{R}(\Gamma,\sigma)$.

\end{proof}

Therefore by the range of the higher trace Theorem 5.3,
and the discussion above, we see that

\begin{thm}[Generalized Quantum Hall Effect] Suppose that the energy level $E$
lies in a gap of the spectrum
of the Hamiltonian $H_\sigma$, then the Hall conductance
$$
\sigma_E = \tr^K (P_E, P_E, P_E) = \tr_c (P_E, P_E, P_E)  \in \phi \mathbb Z
$$
That is, the Hall conductance
has plateaus which are {\em integer} multiples of $\phi$
on any gap in the spectrum of the Hamiltonian, where
$\phi = 2(g-1) + (n-\nu) \in \mathbb Q$.
\end{thm}

\begin{rems}
The set of possible
denominators for $\phi$ is finite
and has been explicitly determined in \cite{Bro}.
It is plausible that this Theorem
might shed light on the mathematical mechanism
responsible for fractional quantum numbers that occur
in the Quantum Hall Effect.

\end{rems}

\end{document}